\documentclass[final,reqno]{siamltex}

\usepackage[english]{babel}
\usepackage{verbatim}
\usepackage{amsmath}
\usepackage{amsfonts}

\usepackage[utf8]{inputenc}
\usepackage{mathrsfs}
\usepackage{paralist}
\usepackage{times}
\usepackage{graphicx}
\usepackage{booktabs}
\usepackage{color}

\usepackage{subfigure}

\newtheorem{example}[theorem]{Example}

\newcommand{\C}{\mathbb{C}}
\newcommand{\R}{\mathbb{R}}

\newcommand{\bI}{\mathbb{I}}

\newcommand{\bM}{\mathbb{M}}
\newcommand{\bN}{\mathbb{N}}

\newcommand{\cK}{\mathcal{K}}
\newcommand{\cL}{\mathcal{L}}

\renewcommand{\i}{\mathrm{i}}

\newcommand{\rand}{\partial}

\newcommand{\abs}[1]{\lvert#1\rvert}
\newcommand{\abslr}[1]{\left\vert #1 \right\vert}

\DeclareMathOperator{\sn}{sn}

\DeclareMathOperator{\re}{Re}
\DeclareMathOperator{\im}{Im}

\newcommand{\coloneq}{\mathrel{\mathop:}=}

\renewcommand{\bI}{{\bf I}}

\renewcommand{\bM}{{\bf M}}
\renewcommand{\bN}{{\bf N}}
\newcommand{\bt}{{\bf t}}

\title{Fast and accurate computation of the logarithmic capacity of compact 
sets}

\author{J\"{o}rg Liesen\thanks{Institute of Mathematics, Technische 
Universit\"{a}t Berlin, Stra{\ss}e des 17. Juni 136, 10623 Berlin, Germany
(\texttt{liesen@math.tu-berlin.de}).}
\and
Olivier S\`{e}te\thanks{Mathematical Institute, University of Oxford, 
Oxford, OX2 6GG, United Kingdom
(\texttt{olivier.sete@maths.ox.ac.uk}).
The second author was supported by the European 
Research Council under the European Union's Seventh Framework Programme 
(FP7/2007-2013)/ERC grant agreement no. 291068.  The views expressed in this 
article are not those of the ERC or the European Commission, and the European 
Union is not liable for any use that may be made of the information contained 
here.}
\and
Mohamed M.S. Nasser\thanks{Department of Mathematics, Statistics and Physics, 
Qatar University, P.O.Box 2713, Doha, Qatar
(\texttt{mms.nasser@qu.edu.qa}).}
}

\begin{document}
\maketitle

\pagestyle{myheadings}
\thispagestyle{plain}
\markboth{J. Liesen, O. S\`{e}te, and M.M.S. Nasser}{Fast and accurate 
computation of the logarithmic capacity of compact sets}

\begin{abstract}
We present a numerical method for computing the logarithmic capacity
of compact subsets of $\C$, which are bounded by Jordan curves and have 
finitely connected complement.
The subsets may have several components and need not have any special symmetry.
The method relies on the conformal map onto lemniscatic domains and, 
computationally, on the solution of a boundary integral equation with 
the Neumann kernel.
Our numerical examples indicate that the method is fast and accurate.
We apply it to give an estimate of the logarithmic capacity of the Cantor 
middle third set and generalizations of it.
\end{abstract}



\section{Introduction}

The \emph{logarithmic capacity} $c(E)$ of a compact set $E$ in the complex plane
${\mathbb C}$ is an important invariant that is closely related to polynomial 
approximation, potential theory, and conformal mapping.  If $E$ is simply 
connected, then $c(E)$ can be obtained via the Riemann map from the complement 
of $E$ 
onto the complement of the unit disk.
The Riemann map is known analytically for some simple sets such as disks, 
ellipses and intervals.  This leads to analytic formulas for the logarithmic 
capacity in terms of the parameters describing the corresponding sets.
Some examples are shown in Table~\ref{table:known_cap} below;
see also~\cite[p.~135]{Ransford1995},~\cite[p.~557]{Ransford2010}
or~\cite[pp.~172--173]{Landkof1972}.

The logarithmic capacity can be obtained from the (exterior) 
Schwarz--Christoffel map when $E$ consists of one or several components with 
polygonal boundaries and has a connected complement; see, 
e.g.,~\cite[Sections~4.4 and~5.8]{DriscollTrefethen2002}.
For simply connected sets $E$ this has been implemented in the 
Schwarz--Christoffel Toolbox~\cite{SCtoolbox}, where the logarithmic capacity 
is one of the outputs of the command \texttt{extermap}.
We refer to~\cite{EmbTre99} for some interesting examples of polygonal sets
with special symmetry properties with respect to the real line.

Considering more complicated sets, in particular the Cantor middle third set, 
Ransford and Rostand pointed out that computing the capacity is 
``notoriously hard''~\cite[p.~1499]{RansfordRostand2007}.
They derived a method for computing upper and lower bounds
for the logarithmic capacity, which can in principle be made arbitrarily close
to each other.  Methods for a direct computation of the logarithmic capacity
have been proposed in~\cite{Dijkstra2009,Rostand1997}; also see the survey
given in~\cite{Ransford2010}. Interest in numerically computing the
logarithmic capacity goes back at least to~\cite{Davis1957}.

In this paper we show that the logarithmic capacity of a fairly wide class
of sets can be computed fast and accurately using conformal mapping techniques. 
The method we present is based on a conformal map from the complement of $E$ 
onto a lemniscatic domain, which is originally due to Walsh~\cite{Walsh1956} 
and which represents a direct generalization of the Riemann map to sets $E$ 
with several components.  We derived a numerical method for computing Walsh's 
map in~\cite{NLS2016_numconf}.  The logarithmic capacity is obtained, 
almost as a by-product, in the first step of that method, and without computing 
the conformal map itself. Because of the practical importance of the logarithmic 
capacity and the apparent lack of general purpose software for its computation, 
we here derive a stand-alone method and its MATLAB implementation.

Our method is applicable to any compact set $E$ whose complement in the extended
plane is connected and bounded by finitely many (sufficiently smooth) Jordan
curves. In particular, it is required neither that $E$ is connected nor that
$E$ has special symmetry properties. Going beyond the theory presented 
in~\cite{NLS2016_numconf}, we place a particular emphasis on the treatment 
of sets $E$ with corners. Moreover, we use a recently developed iterative 
method for mapping parallel slit domains onto domains exterior to ellipses
(see Appendix~A) in order to compute approximations of the logarithmic capacity 
of the classical Cantor middle third set and generalizations of it. Numerous
further numerical examples demonstrate that our method is fast and accurate.
In our numerical examples with sets for which the logarithmic capacity is
known analytically, our method typically yields a computed approximation with a
relative error of order $10^{-14}$.  These computations in MATLAB take at most
a few seconds on a standard laptop.

The paper is organized as follows. In Section~\ref{sect:background} we summarize
major facts about the logarithmic capacity and its relation to conformal mapping.
In Section~\ref{sect:method} we describe our method for computing the logarithmic 
capacity and state its MATLAB implementation. In Section~\ref{sect:numerics} we 
give numerical examples.

\section{Background}
\label{sect:background}

The \textit{transfinite diameter} of a compact set $E \subset \C$ is defined as
\begin{equation}\label{eqn:td}
d(E) \coloneq \lim_{n\rightarrow\infty}
d_n(E)^{\frac{2}{n(n-1)}},\quad\mbox{where}\quad
d_n(E) \coloneq \max_{z_1,\dots,z_n\in E}\prod_{1\leq k<\ell\leq n}
|z_k-z_\ell|.
\end{equation}
The existence of the limit in \eqref{eqn:td} was first shown by
Fekete~\cite{Fek23},
who considered the $d_n$ a sequence of ``generalized diameters'' of $E$.
Note that $d_2(E)=\max_{z_1,z_2\in E}|z_1-z_2|$ is the (usual) diameter of $E$.
In the same article Fekete showed that the transfinite diameter is equal
to the \textit{Chebyshev constant}
\begin{equation}\label{eqn:cc}
t(E) \coloneq \lim_{n\rightarrow\infty}t_n^{\frac{1}{n}},\quad\mbox{where}\quad
t_n \coloneq \min_{\deg(p) \leq n-1} \max_{z\in E} |z^n-p(z)|.
\end{equation}
Thus, the geometric constant $d(E)$ is closely related to polynomial approximation
in the complex plane. This relation can be used to easily show that for both the
Chebyshev constant and the transfinite diameter it is sufficient to consider
compact sets $E$ with a connected complement.  Let $E \subset \C$ be a compact 
set and let $E^c \coloneq \widehat{\C} \backslash E$ denote its complement in 
the extended complex plane $\widehat{\C} = \C \cup \{ \infty \}$.
Denote by $\cK$ the component of the complement that
contains the point at infinity (i.e., $\infty \in \cK$) and  define
$\widehat{E} = \widehat{\C} \backslash \cK$.  Intuitively, $\widehat{E}$ is
obtained after ``filling in the holes'' in $E$.
Now the definition of the Chebyshev constant and the maximum modulus principle
imply that $d(E)=t(E)=t(\widehat{E})=d(\widehat{E})$.

Szeg\H{o}~\cite{Sze24} showed that if the complement of a compact set $E$ 
is connected and regular in the sense that it possesses a Green's function 
$g_{E^c}$ with pole at infinity, then the transfinite 
diameter and the Chebyshev constant are equal to the
\textit{logarithmic capacity}\footnote{More precisely, Szeg\H{o} showed that
$d(E)$ is equal to the \textit{Robin constant} $\gamma$, which he defined via
$\lim_{z \to \infty} (\log \abs{z} - g_{E^c}(z)) = 
\log\gamma$. In the modern literature
the definition of the Robin constant usually gives $c(E)=\exp(-\gamma)$.}
\begin{equation}\label{eqn:cap_green}
c(E) \coloneq \lim_{z \to \infty} \exp( \log \abs{z} - g_{E^c}(z) ).
\end{equation}
Saff~\cite{Saf10} called the result $d(E)=t(E)=c(E)$ the
\textit{fundamental theorem of classical potential theory}.

If $E^c = \widehat{\C} \backslash E$ is simply connected (and $E$ is not a 
single point), then its Green's function is given by $g_{E^c}(z) = \log 
\abs{\Phi(z)}$, where 
$w=\Phi(z)$ is the uniquely determined Riemann map
from the complement of $E$ to the complement of the unit disk that is normalized
by $\Phi(\infty)=\infty$ and $\Phi'(\infty)>0$. Near infinity this map can be
written as
\begin{equation*}
\Phi(z)=\frac{z}{\mu}+\mu_0+ O \left(\frac{1}{z}\right),
\end{equation*}
where $1/\Phi'(\infty)=\mu>0$ is called the \textit{conformal radius} of $E$.
Hence we have
\begin{equation*}
d(E)=t(E)=c(E)=\lim_{z \to \infty} \exp(\log \abs{z} -\log\abs{\Phi(z)} )=\mu,
\end{equation*}
which shows that the logarithmic capacity can be obtained using (analytical or
numerical) conformal mapping techniques.  The simplest example is a disk $E$ of
radius $r>0$, for which $\Phi(z)=z/r$ and thus $c(E)=r$; see
Table~\ref{table:known_cap} for some further analytically known examples.

\begin{table}[t]
\begin{center}
\begin{tabular}{lc}
\toprule
$E$ & $c(E)$ \\
\midrule
disk of radius $r$ & $r$ \\[0.2ex]
half-disk of radius $r$ & $4r/3^{3/2}$ \\[0.2ex]
ellipse with semi-axes $a$ and $b$ & $\frac12(a+b)$ \\[0.2ex]
line segment of length $h$ & $\frac14 h$ \\[0.2ex]
square with side $h$ & $ \frac14 \Gamma(1/4)^2h/\pi^{3/2}$ \\[0.2ex]
two intervals $[-b, -a] \cup [a,b]$ & $\frac12 \sqrt{b^2-a^2} $ \\[0.2ex]
\bottomrule
\end{tabular}
\end{center}
\caption{Examples of known logarithmic capacities.}
\label{table:known_cap}
\end{table}

Walsh~\cite{Walsh1956} proved a direct generalization of the classical Riemann
mapping theorem in which he replaced the complement of the unit disk
by a \emph{lemniscatic domain} of the form
\begin{equation}\label{eqn:lemniscate}
\cL \coloneq \{ z \in \widehat{\C} : \abs{U(z)} > \mu \}, \quad \text{where}
\quad U(z) \coloneq \prod_{j=1}^\ell (z-a_j)^{m_j},
\end{equation}
$a_1, \ldots, a_\ell \in \C$ are pairwise distinct, $m_1, \ldots, m_\ell$
are positive real numbers with $\sum_{j=1}^\ell m_j = 1$, and $\mu > 0$.
A simple calculation shows that $\log\abs{U(z)} - \log(\mu)$ is the Green's
function with pole at infinity for $\cL$, so that by~\eqref{eqn:cap_green} we have
\begin{equation*}
c( \widehat{\C} \backslash \cL ) = \lim_{z \to \infty} \exp( \log\abs{z} -
\log\abs{U(z)} + \log(\mu) ) = \mu.
\end{equation*}
Walsh proved the following existence theorem in~\cite{Walsh1956}.

\begin{theorem}\label{thm:Walsh}
Let $E$ be a compact set whose complement $\cK = \widehat{\C} \backslash E$ is
connected and bounded by $\ell$ Jordan curves. Then there exists a uniquely determined
lemniscatic domain $\cL$ of the form~\eqref{eqn:lemniscate} with $\mu=c(E)$
and a uniquely determined conformal map
\begin{equation*}
\Phi : \cK \to \cL \quad \text{with} \quad \Phi(z) = z + O
\left(\frac{1}{z}\right) \text{ near infinity.}
\end{equation*}
\end{theorem}

The first analytic examples of Walsh's conformal map onto lemniscatic domains
have recently been given in~\cite{SeteLiesen2016_conf}.  These were
applied in~\cite{SeteLiesen2016_fwprop} in a study of polynomial approximation
problems on disconnected compact sets and, in particular, two real intervals.
In~\cite{NLS2016_numconf} we developed a numerical method for computing the 
lemniscatic domain $\cL$ and the conformal map $\Phi$ corresponding to a 
given compact set $E$ with connected complement $\cK$.
 
\section{Computing the logarithmic capacity}
\label{sect:method}

As indicated in the formulation of Theorem~\ref{thm:Walsh}, the logarithmic
capacity occurs naturally as one of the parameters defining the lemniscatic
domain $\cL$ onto which $\cK = \widehat{\C} \backslash E$ is mapped.
Hence the map $\Phi$ itself is not required for computing $c(E)$. A closer
inspection of the numerical method developed in~\cite{NLS2016_numconf} shows
that the parameters $m_1,\dots,m_\ell,\mu$ of $\cL$ are indeed computed in a
step that can be executed separately. Due to the structure of the underlying
equations the computation of these parameters can be done in a very efficient
way. We will now briefly describe this computation. A detailed derivation is
given in~\cite{NLS2016_numconf}.

As in Theorem~\ref{thm:Walsh}, let $E$ be compact with a finitely
connected complement $\cK$.
(As noted in Section~\ref{sect:background}, the connectedness of $\cK$ is
no restriction for computing the capacity.)
We assume that the boundary $\Gamma = \partial E = \partial \cK$ of $E$
consists of $\ell$ Jordan curves $\Gamma_1, \ldots, \Gamma_\ell$, which satisfy
the following smoothness assumption:
Each $\Gamma_j$ is parameterized by a $2\pi$-periodic function $\eta_j : J_j
\coloneq [0, 2\pi] \to \Gamma_j$, which is twice continuously
differentiable and satisfies $\dot{\eta}_j(t) = \frac{d \eta_j}{dt}(t) \neq 0$
for all $t$.
These assumptions can be relaxed so as to include domains with corners, see
the precise statement and discussion below.
The boundary of $E$ is oriented clockwise, so that $\cK$ is to
the left of the boundary.

Then a parameterization for the whole boundary $\Gamma$ is given by the map
\begin{equation}
\eta : J \to \Gamma = \bigcup_{j=1}^\ell \Gamma_j, \quad
\eta(t) = \begin{cases} \eta_1(t), & t \in J_1, \\ \quad \vdots & \\
\eta_\ell(t), &
t \in J_\ell, \end{cases}
\label{eqn:eta}
\end{equation}
where $J$ is the disjoint union of the intervals $J_1, \ldots, J_\ell$, i.e.,
$J$ consists of $\ell$ copies of $[0, 2\pi]$.

For each $j = 1, 2, \ldots, \ell$, we choose an auxiliary point $\alpha_j$ in
the interior of the Jordan curve $\Gamma_j$, and define the function
\begin{equation}
\gamma_j(t) = - \log \abs{\eta(t)-\alpha_j}, \quad t \in J.
\label{eqn:gammaj}
\end{equation}
In practical applications of our method the parameters $\alpha_j$ must be
specified by the user; cf. the MATLAB code shown in Figure~\ref{fig:logcapacity}
below. Our numerical experience with the method suggests that the actual values
of the $\alpha_j$ are not important, as long as these points are sufficiently
far away from the boundary $\Gamma_j$. In the experiments discussed in
Section~\ref{sect:numerics} we always chose $\alpha_j$ close to (or at)
the center of the interior of $\Gamma_j$. The (blue) dots in
Figures~\ref{fig:bw_hd} and~\ref{fig:multcomp} show some examples.

Let $H$ denote the space of all functions $f$ in $J$, whose restriction to
$J_j = [0, 2\pi]$ is a real-valued, $2\pi$-periodic and H\"{o}lder continuous
function for each $j = 1, 2, \ldots, \ell$.
Define the integral operators
\begin{align*}
(\bN f)(s) &= \int_J \frac{1}{\pi} \im \left(
\frac{\dot{\eta}(t)}{\eta(t)-\eta(s)} \right) f(t) \, dt, \quad s \in J, \\
(\bM f)(s) &= \int_J \frac{1}{\pi} \re \left(
\frac{\dot{\eta}(t)}{\eta(t)-\eta(s)} \right) f(t) \, dt, \quad s \in J,
\end{align*}
on $H$.  The kernel of $\bN$ is called the \emph{Neumann kernel}.  Denoting by
$\bI$ the identity operator on $H$, we can state the following
theorem that combines~\cite[Theorems~4.1 and~4.3]{NLS2016_numconf}.

\begin{theorem}
\label{thm:linsys_cap}
For each $j = 1, 2, \ldots, \ell$ the integral equation
\begin{equation}
(\bI - \bN) \mu_j = - \bM \gamma_j \label{eqn:bieq}
\end{equation}
with $\gamma_j$ as in \eqref{eqn:gammaj} has a unique solution $\mu_j \in H$,
and the function
\begin{equation}
h_j \coloneq ( \bM \mu_j - (\bI - \bN) \gamma_j )/2 \label{eqn:def_hj}
\end{equation}
is real-valued and piecewise constant, that is
\begin{equation*}
h_j(t) = h_{k,j}, \quad t \in J_k, \quad k = 1, 2, \ldots, \ell.
\end{equation*}
Furthermore, $\log (\mu)$ and the parameters $m_1, \ldots, m_\ell$
in~\eqref{eqn:lemniscate} are the unique solution of the linear
algebraic system
\begin{equation}
\begin{bmatrix}
h_{1,1} & h_{1,2} & \cdots & h_{1,\ell} & -1 \\
h_{2,1} & h_{2,2} & \cdots & h_{2,\ell} & -1 \\
\vdots & \vdots & \ddots & \vdots & \vdots \\
h_{\ell,1} & h_{\ell,2} & \cdots & h_{\ell,\ell} & -1 \\
1 & 1 & \cdots & 1 & 0 \\
\end{bmatrix}
\begin{bmatrix} m_1 \\ m_2 \\ \vdots \\ m_\ell \\ \log(\mu) \end{bmatrix}
= \begin{bmatrix} 0 \\ 0 \\ \vdots \\ 0 \\ 1 \end{bmatrix}. \label{eqn:linsys}
\end{equation}
\end{theorem}

This suggests the following method for computing the logarithmic
capacity $\mu$:

\begin{compactenum}
\renewcommand{\labelenumi}{(\arabic{enumi})}
\item \label{it:one}
For each $j=1,\dots,\ell$ solve the integral equation \eqref{eqn:bieq}
for the unknown function $\mu_j$.
\item \label{it:two}
Solve the linear algebraic system \eqref{eqn:linsys} of order $\ell+1$,
where the entries $h_{k,j}$ in the coefficient matrix are computed from
\eqref{eqn:def_hj} using the known functions $\gamma_j$ and the functions
$\mu_j$ computed in step (\ref{it:one}).
\end{compactenum}

The linear algebraic system in step (\ref{it:two}) is usually quite small and
we solve it directly using the ``backslash'' operator in MATLAB. Step
(\ref{it:one}) requires more work:

As described in~\cite[Section~5]{NLS2016_numconf}, the $\ell$ boundary
integral equations \eqref{eqn:bieq} can be solved accurately by the
Nystr\"om method with the trapezoidal rule. This method
yields a linear algebraic system with a dense nonsymmetric matrix $I-B$
of order $\ell n$, where $n$ is the number of nodes in the discretization
of each boundary component. This system can be solved iteratively
using the GMRES method. Each step of this method requires one multiplication
with the matrix $I-B$. Due to the structure of the integral equation, this
product can be computed efficiently in $O(\ell n)$ operations using the Fast
Multipole Method (FMM). The eigenvalues of the matrix $I-B$ are contained
in the interval $(0,2]$ and they cluster around $1$. As observed
in~\cite{NLS2016_numconf} and several other publications (see,
e.g.,~\cite{Nas-fast,Nas-siam13}), the number of GMRES iterations
for obtaining a very good approximation of the exact solution is
mostly independent of the given domain and number of nodes in the
discretization of the boundary.

The method for solving \eqref{eqn:bieq} for the $\mu_j$ and
subsequently computing $h_j$ in \eqref{eqn:def_hj} has been implemented
in the MATLAB function {\tt fbie} shown in~\cite[Fig.~4.1]{Nas-fast}.
This function uses MATLAB's built-in {\tt gmres} function as well
as the function {\tt zfmm2dpart} from the fast multipole toolbox
FMMLIB2D~\cite{Gre-Gim12}. The main inputs of the method consist of
the discretized functions $\eta(t)$, $\dot{\eta}(t)$, and $\gamma_j(t)$
described above.

For domains with smooth boundaries, we discretize the interval
$[0, 2\pi]$ by $n$ nodes $s_1, \ldots, s_n$ and write $s = [s_1, \ldots, s_n]$
where $n$ is an even integer.  For simplicity, we
usually take equidistant nodes, i.e.,
\begin{equation}\label{eq:s_i}
s_k = (k-1) \frac{2 \pi}{n}, \quad k = 1, \ldots, n.
\end{equation}
Then $\ell$ copies of $s$ give a discretization
\begin{equation}\label{eq:bt}
\bt = [s, s, \ldots, s]^T \in \C^{\ell n}
\end{equation}
of the parameter domain $J$, leading to the discretizations
\begin{equation}\label{eq:dis_s}
\eta(\bt) = [\eta_1(s), \eta_2(s), \ldots, \eta_\ell(s)]^T, \quad
\dot{\eta}(\bt), \quad
\gamma_j(\bt) = - \log \abs{ \eta(\bt) - \alpha_j } \in \C^{\ell n},
\end{equation}
$j=1,2,\ldots,\ell$.
We store the discretized functions in the vectors \texttt{et}, \texttt{etp},
\texttt{gam}, respectively, and call
\begin{verbatim}
  [~,h] = fbie(et,etp,ones(size(et)),gam,n,5,[],tol,maxit)
\end{verbatim}
Here {\tt []} means that GMRES is used without restart,
{\tt tol} is the convergence criterion used within GMRES, and {\tt maxit}
is the maximal number of GMRES iterations. In the numerical experiments
described in Section~\ref{sect:numerics} we have used {\tt tol=1e-14} and
{\tt maxit=100}.  The output of \texttt{fbie} are the values $h_j(\bt)$ of
$h_j$ from~\eqref{eqn:def_hj}, and the values $h_{k,j}$ are computed
by taking arithmetic means:
\begin{equation*}
h_{k,j} = \frac{1}{n} \sum_{i=1+(k-1)n}^{kn} h_j(t_i), \quad k = 1, 2, \ldots,
\ell, \quad j = 1, 2, \ldots, \ell.
\end{equation*}
These values are used to set up the linear algebraic system~\eqref{eqn:linsys}, 
which we solve directly as mentioned above.

Figure~\ref{fig:logcapacity} shows our MATLAB implementation of the overall
method described in this section, where the inputs \texttt{et} and \texttt{etp}
are given as in~\eqref{eq:dis_s}, and \texttt{alpha} is the column vector of
auxiliary points $\alpha_j$.

\begin{figure}
\begin{center}
\footnotesize
\begin{verbatim}
function mu = logcapacity(et,etp,alpha)
% Computes the logarithmic capacity of a compact set with L components.
%
% Input:
%  et    = [eta_1; ...; eta_L] (discretized boundary; column vector)
%  etp   = derivative of the boundary curves (same format as et)
%  alpha = [alpha(1); ...; alpha(L)] (auxiliary point alpha(j) interior
%                                     to j-th boundary curve Gamma_j)

L = length(alpha);  %% number of boundary components
n = length(et)/L;   %% number of nodes per boundary component

% Auxiliary functions gamma_j(t) = -log |eta(t)-alpha_j|
for k=1:L
    gamj(:,k) = -log(abs(et-alpha(k)));
end

% Compute the auxiliary functions h_j
A = ones(size(et)); %% the function A in the gen. Neumann kernel
for k=1:L
    [~,hjv(:,k)] = fbie(et,etp,A,gamj(:,k),n,5,[],1e-14,100);
end

% Build and solve linear system for m_1, ..., m_L, log(mu)
for j=1:L
    for k=1:L
        hj(k,j) = sum(hjv(1+(k-1)*n:k*n,j))/n;
    end
end
matA            = hj;
matA(L+1,1:L)   = 1;
matA(1:L,L+1)   = -1;
vc_right        = zeros(L+1,1);
vc_right(L+1)   = 1;

x  = matA\vc_right;
mu = exp(x(L+1));
end
\end{verbatim}
\end{center}
\caption{MATLAB code for the computation of the logarithmic capacity.}
\label{fig:logcapacity}
\end{figure}

The presented method can be extended to domains with corners.
We assume that the corner points are not cusps and that the tangent vector of
the boundary has only the first kind discontinuity at these corner points.
The left tangent vector at a corner point is considered as the tangent vector
at this point.  In this case, the solution of the integral
equation~\eqref{eqn:bieq} has a singularity in its first derivative in the
vicinity of the corner points~\cite{Nas-corner}.  Using the equidistant
nodes~\eqref{eq:bt} to discretize the integrals in~\eqref{eqn:bieq}
and~\eqref{eqn:def_hj} yields
only poor convergence~\cite{Kress1990,Nas-corner,Rathsfeld93}.
To achieve a satisfactory accuracy, we discretize
the integrals in~\eqref{eqn:bieq} and~\eqref{eqn:def_hj} using a \emph{graded
mesh} which is based on substituting a new variable in such a way that
the discontinuity of the derivatives of the solution of the integral
equation at the corner points is removed.

Following Kress~\cite{Kress1990}, we define a bijective, strictly monotonically
increasing and infinitely differentiable function,
$w : [0,2\pi] \to [0,2\pi]$, by
\begin{equation*}
w(t)=2\pi\frac{v(t)^p}{v(t)^p+v(2\pi-t)^p},
\end{equation*}
where
\begin{equation*}
v(t)=\left(\frac{1}{p}-\frac{1}{2}\right)\left(\frac{\pi-t}{\pi}\right)^3+\frac{
1}{p}\frac{t-\pi}{\pi}+\frac{1}{2}.
\end{equation*}
The grading parameter is the integer $p\ge2$, and the cubic
polynomial $v$ is chosen to ensure that, for the equidistant mesh
$s_1,s_2,\ldots,s_n$, almost $n/2$ of the grid points
$w(s_1), w(s_2), \ldots$, $w(s_n)$ are equally distributed throughout
$[0,2\pi]$, and the other half is accumulated towards the two endpoints $0$ and
$2\pi$.

Assume that the boundary $\Gamma_k$ has $q_k>0$ corner points
\begin{equation*}
\eta_k(0),\quad \eta_k(2\pi/q_k),\quad \eta_k(4\pi/q_k),
\quad\ldots, \quad \eta_k(2(q_k-1)\pi/q_k).
\end{equation*}
Then we define a bijective, strictly monotonically increasing and infinitely
differentiable function, $\delta_k:J_k\to J_k$, by
\[
\delta_k(t)= \begin{cases}
 w(q_k t)/q_k,                      & t \in [0,2\pi/q_k), \\
(w(q_k t-2\pi)+2\pi)/q_k,    & t \in [2\pi/q_k,4\pi/q_k), \\
(w(q_k t-4\pi)+4\pi)/q_k,    & t \in [4\pi/q_k,6\pi/q_k), \\
\quad \vdots & \\
(w(q_k t-2(q_k-1)\pi)+2(q_k-1)\pi)/q_k,
& t \in [2(q_k-1)\pi/q_k,2\pi].
\end{cases}
\]
Since the function $w$ has a zero of order $p$ at the endpoints $t=0$ and
$t=2\pi$~\cite[Theorem~2.1]{Kress1990}, the function $\delta_k$ is at least $p$
times continuously differentiable.
If the boundary $\Gamma_k$ has no corner points we define the function
$\delta_k(t)$ by
\[
\delta_k(t)=t, \quad t\in J_k.
\]
Finally, we define a function $\delta:J\to J$ by
\begin{equation*}
\delta(t) =
\begin{cases}
\delta_1(t),    & t \in J_1, \\
\quad \vdots    & \\
\delta_\ell(t), & t \in J_\ell.
\end{cases}
\end{equation*}

As noted in the first paragraph in~\cite[p.~242]{Kress1991}, using
the graded mesh
\begin{equation}\label{eq:bt-cr}
\delta(\bt) = [\delta_1(s), \delta_2(s), \ldots, \delta_\ell(s)]^T \in \C^{\ell
n}
\end{equation}
for discretizing the integrals
in~\eqref{eqn:bieq} and~\eqref{eqn:def_hj} is equivalent to parameterizing
the boundary $\Gamma$ by $\eta(\delta(t))$, and then solving the integral
equation
as in the case of smooth domains. Hence, we have the discretizations
\begin{equation}\label{eq:dis_s-cr}
\eta(\delta(\bt)), \quad
\dot{\eta}(\delta(\bt)) \dot{\delta}(\bt), \quad
\gamma_j(\delta(\bt)) = - \log \abs{ \eta(\delta(\bt)) - \alpha_j } \in \C^{\ell
n},
\end{equation}
$j=1,2,\ldots,\ell$, which replace the discretized functions
in~\eqref{eq:dis_s}.
In our numerical experiments we have used the grading parameter $p=3$.

 
\section{Numerical examples}
\label{sect:numerics}

We now present numerical examples that illustrate our meth\-od.  If
not stated otherwise, we use an equidistant discretization of $[0, 2\pi]$.
Computations were performed in MATLAB R2013a on an ASUS Laptop with Intel Core
i7-4720HQ CPU @ 2.60Ghz 2.59 Ghz and 16~GB RAM using the code shown
in Figure~\ref{fig:logcapacity}.  Computation times were measured with the
MATLAB \verb|tic toc| command.

\subsection{Sets with one component}

\begin{example}[Disk] \label{ex:disk}
{\rm
Let $E_r \coloneq \{ z \in \C : \abs{z} \leq r \}$ denote the closed disk with
radius
$r$ and $c(E_r)=r$; see Table~\ref{table:known_cap}. We use the parametrization
\begin{equation*}
\eta : [0, 2\pi] \to \rand E_r, \quad t \mapsto r e^{-\i t},
\end{equation*}
and $n = 2^8 = 256$ points in the discretization. For $r = 1$ our method
computes the exact value $c(E_1) = 1.0$. For $r = 2$ it computes the value
$2.000000000000003$, accurate to $14$ digits.  (The relative error is $1.33
\cdot 10^{-15}$.)  Each computation took less than $0.1~\mathrm{s}$.
}
\end{example}

\begin{example}[Ellipse] \label{ex:ellipse}
{\rm
We consider the family of ellipses $E_d$ with semi-axes $a=1$ and
$b=10^{-d}$, and $c(E_d) = (1+10^{-d})/2$; see Table~\ref{table:known_cap}.
A parametrization of the boundary is
\begin{equation*}
\eta_d(t) = \cos(t) - \i 10^{-d} \sin(t), \quad t \in [0, 2\pi].
\end{equation*}
We use $d = 1, 2, 3, 4$ and $n = 2^k$ with $k = 8, 9, \ldots, 18$.
For $d = 1$ we have $c(E_1) = 0.55$, and the value computed by our method
has a relative error smaller than $10^{-13}$ for every $n$.
For $d = 2, 3, 4$ the relative errors of the computed values are shown
in Figure~\ref{fig:ellipse} (left). We observe that our method is less
accurate for larger $d$, i.e., for ``thinner'' ellipses. For large $d$
the method still yields a very accurate approximation of $c(E_d)$, but
this requires a large increase of the number of discretization
points. This may be related to the fact that the auxiliary point inside
a ``thinner'' ellipse is necessarily closer to the boundary of the ellipse
(cf. our comments after equation \eqref{eqn:gammaj}). On the right of
Figure~\ref{fig:ellipse} we show the computation times.
\begin{figure}[t]
\centerline{
\includegraphics[width=0.5\textwidth]{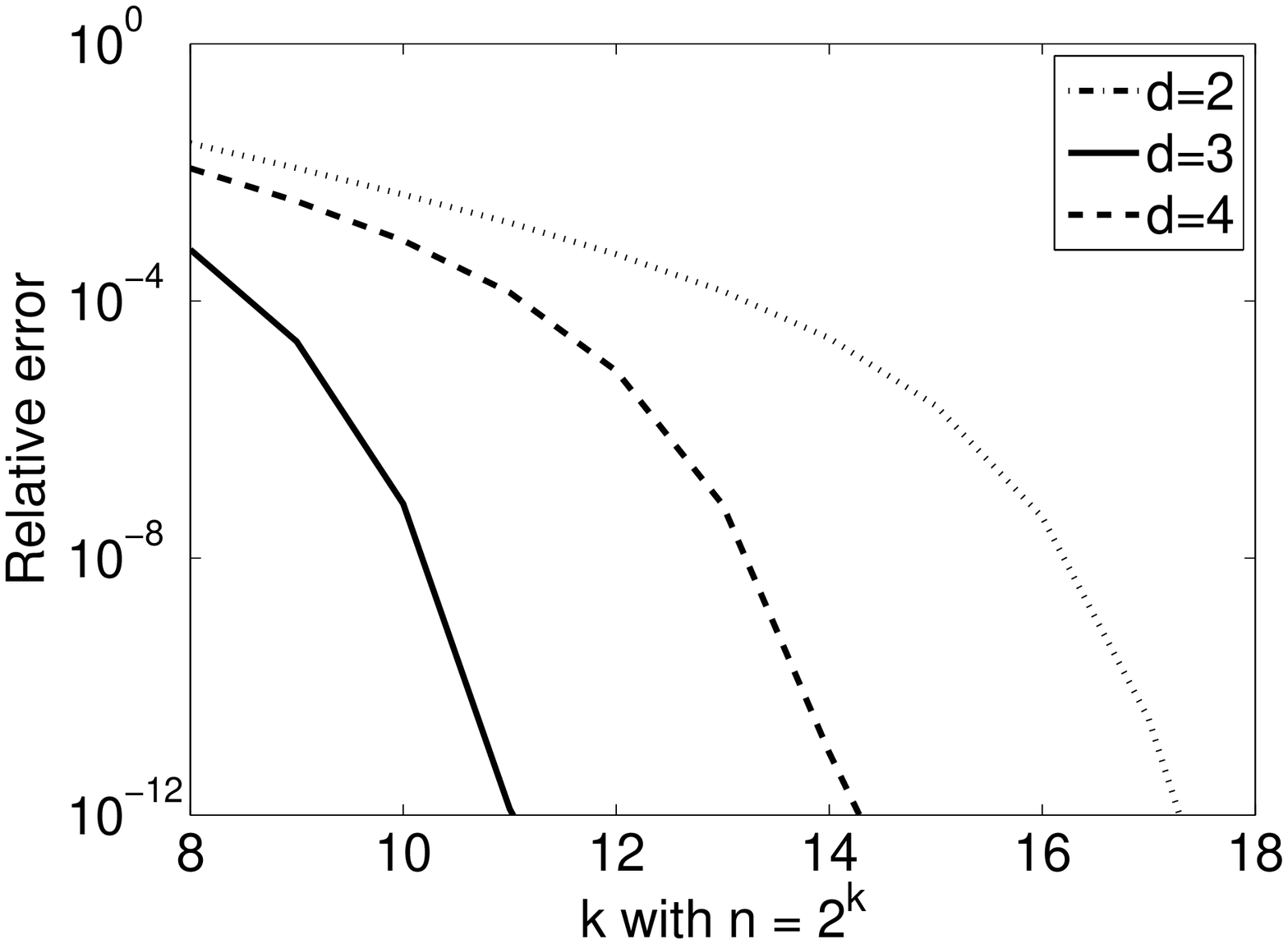}
\includegraphics[width=0.5\textwidth]{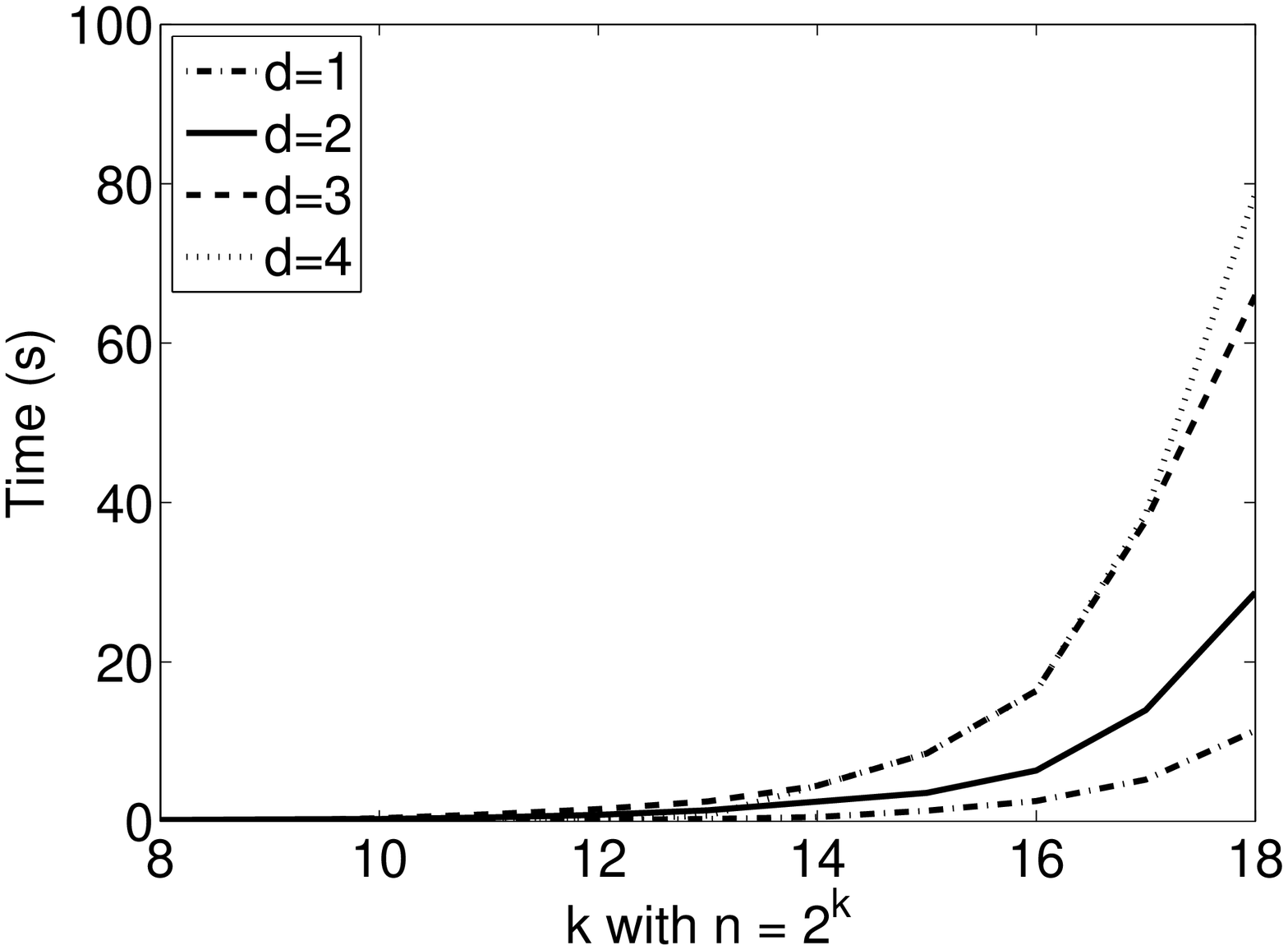}
}
\caption{Results for Example~\ref{ex:ellipse} (Ellipse): Relative errors
of the computed logarithmic capacity (left) and computation times
(in seconds, right).}
\label{fig:ellipse}
\end{figure}
}
\end{example}

\begin{example}[Half-disk] \label{ex:half-disk}
{\rm
Let $E \coloneq \{ z \in \C : \abs{z} \leq 1, \im(z) \geq 0 \}$ be the upper
half
of the unit disk with $c(E) = 4/3^{3/2}$; see Table~\ref{table:known_cap}.
Since the boundary of $E$ has corners, we use a graded mesh as described in
Section~\ref{sect:method}. The following table shows the computed values
(correct digits are \underline{underlined}), relative errors, and computation
times for increasing $n$:
\begin{center}
\begin{tabular}{lllc}
\toprule
$n$ & computed capacity & relative error & time (s) \\
\midrule
$2^{10}$ & $\underline{0.7698003}47826294$ & $1.44 \cdot 10^{-8}$  & $0.1$ \\
$2^{11}$ & $\underline{0.76980035}7536609$ & $1.80 \cdot 10^{-9} $ & $0.2$ \\
$2^{12}$ & $\underline{0.769800358}746887$ & $2.24 \cdot 10^{-10}$ & $0.4$ \\
$2^{13}$ & $\underline{0.769800358}897939$ & $2.80 \cdot 10^{-11}$ & $0.5$ \\
$2^{14}$ & $\underline{0.76980035891}6802$ & $3.51 \cdot 10^{-12}$ & $0.8$ \\
$2^{15}$ & $\underline{0.769800358919}097$ & $5.25 \cdot 10^{-13}$ & $1.5$ \\
$2^{16}$ & $\underline{0.769800358919}667$ & $2.15 \cdot 10^{-13}$ & $3.7$ \\
$2^{17}$ & $\underline{0.7698003589195}14$ & $1.63 \cdot 10^{-14}$ & $6.1$ \\
\bottomrule
\end{tabular}
\medskip
\end{center}
}
\end{example}

\begin{example}[Square] \label{ex:square}
{\rm
Let $E$ be the square with vertices $1+\i, 1-\i, -1-\i$, $-1-\i$ and
hence side length $h=2$, giving $c(E)=\frac12 \Gamma(1/4)^2/\pi^{3/2}$;
see Table~\ref{table:known_cap}. Evaluating this expression in MATLAB
gives $c(E)= 1.180340599016096$, which we use as the ``exact'' value for
our experiment. As for the half-disk, our discretization uses a graded mesh.
The following table shows the computed values (correct digits are
\underline{underlined}), relative errors and computation times for
increasing $n$:
\begin{center}
\begin{tabular}{lllc}
\toprule
$n$ & computed capacity & relative error & time (s) \\
\midrule
$2^8$ & $\underline{1.1803}28330582103$ & $1.04 \cdot 10^{-05}$    & $0.1$ \\
$2^9$ & $\underline{1.1803}39089365394$ & $1.28 \cdot 10^{-06}$    & $0.1$ \\
$2^{10}$ & $\underline{1.180340}411967884$ & $1.58 \cdot 10^{-07}$ & $0.1$ \\
$2^{11}$ & $\underline{1.1803405}75744745$ & $1.97 \cdot 10^{-08}$ & $0.2$ \\
$2^{12}$ & $\underline{1.18034059}6114299$ & $2.46 \cdot 10^{-09}$ & $0.2$ \\
$2^{13}$ & $\underline{1.18034059}8653831$ & $3.07 \cdot 10^{-10}$ & $0.5$ \\
$2^{14}$ & $\underline{1.18034059}8970863$ & $3.83 \cdot 10^{-11}$ & $0.6$ \\
$2^{15}$ & $\underline{1.18034059901}0508$ & $4.73 \cdot 10^{-12}$ & $1.1$ \\
$2^{16}$ & $\underline{1.18034059901}5215$ & $7.47 \cdot 10^{-13}$ & $2.3$ \\
$2^{17}$ & $\underline{1.180340599016}100$ & $3.57 \cdot 10^{-15}$ & $5.9$ \\
\bottomrule
\end{tabular}
\medskip
\end{center}
We also compute the logarithmic capacity with the Schwarz--Christoffel
Toolbox~\cite{SCtoolbox}. With the default settings the commands
\begin{verbatim}
  p = polygon([1+i,-1+i,-1-i,1-i]);
  f = extermap(p);
  capacity(f)
\end{verbatim}
yield the value $\underline{1.1803405990}90706$, which has the relative error
$6.32 \cdot 10^{-11}$.  For a tolerance of $10^{-14}$ instead
of the default value $10^{-8}$ the Schwarz--Christoffel Toolbox returns
a very accurate result with the relative error $3.76 \cdot 10^{-16}$.
}
\end{example}

In both Example~\ref{ex:half-disk} and~\ref{ex:square} the relative error in 
our method converges to the machine precision for increasing $n$.  The reason 
that we need many points to obtain very accurate results is that in both 
examples the boundaries have corners.

\begin{figure}
\centerline{
\subfigure{\includegraphics[width=0.5\textwidth]{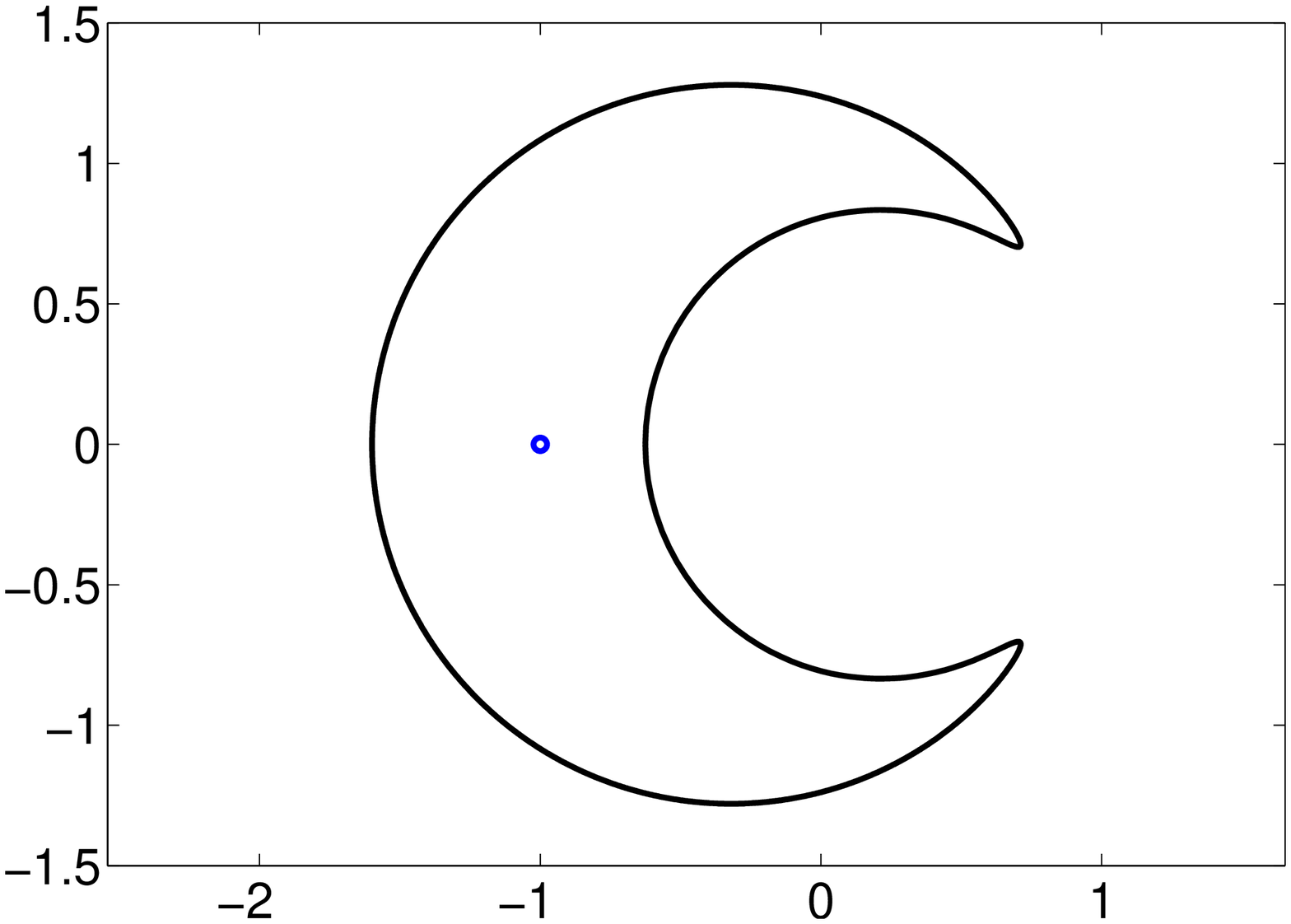}}
\subfigure{\includegraphics[width=0.5\textwidth]{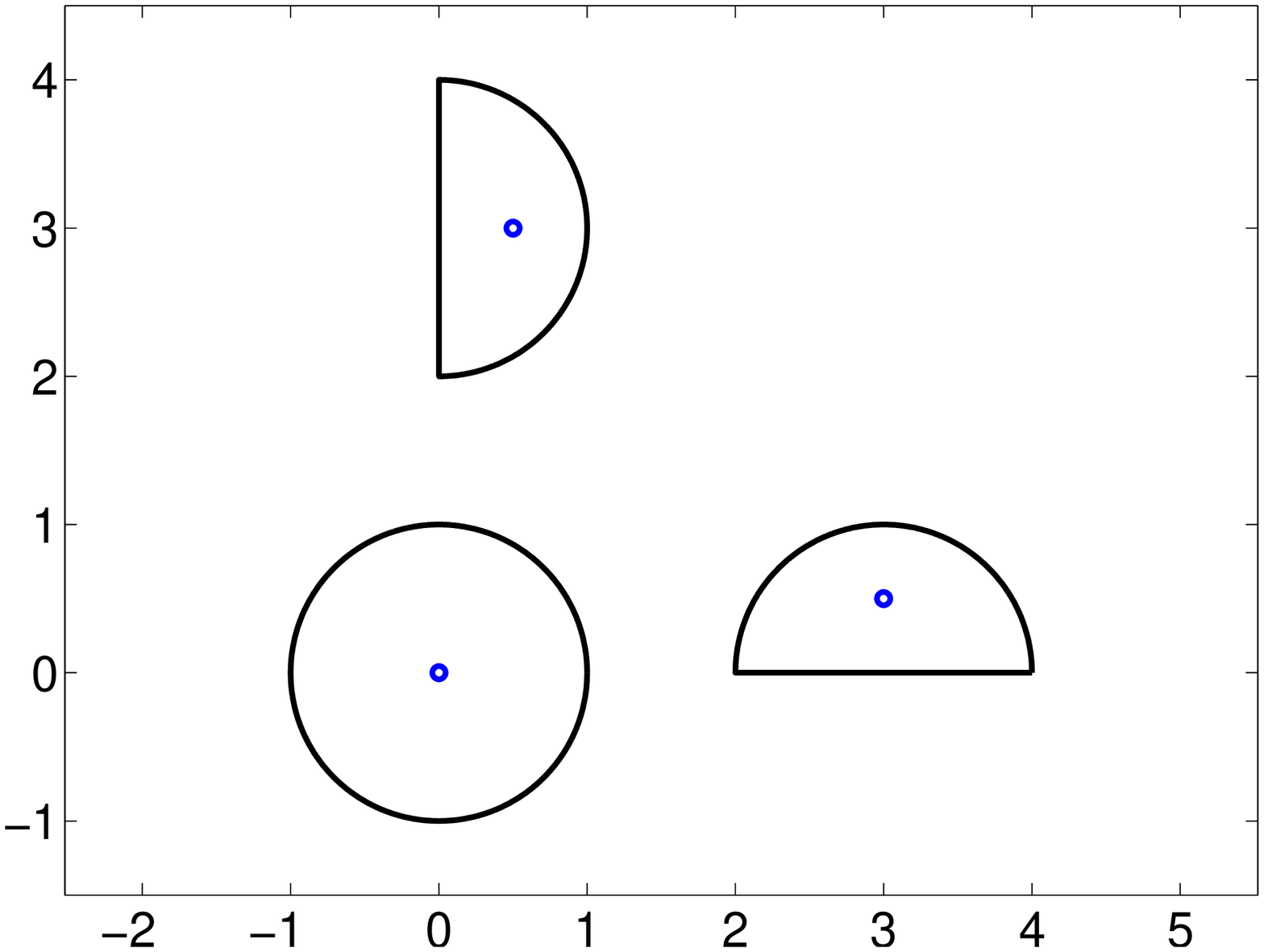}}
}
\caption{The sets from Examples~\ref{ex:bw} (left) and~\ref{ex:Rostand_disks}
(right). The (blue) dots show the auxiliary points $\alpha_j$.}
\label{fig:bw_hd}
\end{figure}

\begin{example} \label{ex:bw}
{\rm
We consider the set in Figure~\ref{fig:bw_hd} (left), which is of the form
introduced in~\cite{KochLiesen2000}.  For
these sets an analytic parameterization of the boundary and the logarithmic
capacity are known explicitly. Here we consider the compact set $E$ bounded by
\begin{equation*}
\eta(t) = \psi(e^{-\i t}), \quad t \in [0, 2\pi],
\end{equation*}
where $\psi$ is the conformal map given by~\cite[Equation~(3.2)]{KochLiesen2000}
with the parameters $\lambda_m = 1$, $\phi = \pi/2$ and $\epsilon = 0.1$.
Then $c(E) = 1.223502096192244$; see~\cite[Corollary~3.4]{KochLiesen2000}.
The following table shows the computed values (correct digits are
\underline{underlined}), relative errors and computation times
for increasing $n$:
\begin{center}
\begin{tabular}{lccc}
\toprule
$n$ & computed capacity & relative error & time (s) \\
\midrule
$2^8$ & $\underline{1.223}385602611070$ & $9.52 \cdot 10^{-05}$ & $0.1$ \\
$2^9$ & $\underline{1.22350}0703601890$ & $1.14 \cdot 10^{-06}$ & $0.1$ \\
$2^{10}$ & $\underline{1.22350209}5786708$ & $3.31 \cdot 10^{-10}$ & $0.1$ \\
$2^{11}$ & $\underline{1.22350209619224}5$ & $3.63 \cdot 10^{-16}$ & $0.2$ \\
\bottomrule
\end{tabular}
\end{center}
As in the two previous examples, the relative error converges to machine
precision as $n$ increases.  The reason that we need many points to obtain an
accurate result is different from the two previous examples.  Here the boundary
is analytic, but with equispaced points in $[0, 2\pi]$ and the above
parametrization, only few discretization points lie on the inner arc, which is
not well resolved for smaller $n$.  Here, a different parametrization might
lead to very accurate results already for small $n$, but we did not pursue
this further.
}
\end{example}

\subsection{Sets with several components}

\begin{example}[Two disks with equal radii] \label{ex:2disks}
{\rm
Let $r, z_0 \in \R$ with $0 < r < z_0$ and let $E$ be the union of the two
disks $D_r(z_0) = \{ z \in \C : \abs{z-z_0} \leq r \}$ and
$D_r(-z_0) = \{ z \in \C : \abs{z+z_0} \leq r \}$.
From~\cite[Theorem~4.2]{SeteLiesen2016_conf} we know that
\begin{equation*}
c(E) = \frac{2 K}{\pi} \sqrt{z_0^2 - r^2} \sqrt{2 L (1 + L^2)},
\end{equation*}
where
\begin{equation*}
\rho = \frac{\sqrt{z_0 + r} - \sqrt{z_0 - r}}{\sqrt{z_0 + r} + \sqrt{z_0 - r}},
\quad
L = 2 \rho \prod_{k=1}^\infty \left( \frac{1 + \rho^{8k}}{1 + \rho^{8k-4}}
\right)^2,
\end{equation*}
and
\begin{equation*}
K = K(L^2) = \int_0^1 \frac{1}{\sqrt{(1-t^2) (1-L^4 t^2)}} \, dt.
\end{equation*}
The product for $L$ converges very quickly, so that it suffices to compute the
first few factors to obtain the correct value up to the machine precision.  Using
this value of $L$ we evaluate the complete elliptic integral of the first kind
for $K$ with the MATLAB command \verb|ellipk|.  We use the result as the
``exact'' value $c(E)$.
The constant $L$ can also be written with the Jacobi theta functions as
$L = \theta_2(0; \rho^4) / \theta_3(0; \rho^4)$, using their product
representation~\cite[§21.3]{WhittakerWatson1962}.
We apply our algorithm for $z_0 = 1$ and $r = 0.5, 0.7,0.9$ with $n = 2^8$
(giving $2n = 2^9$ nodes in total).  The results are very accurate ($14$
digits of accuracy) and are shown in
the following table (correct digits are
\underline{underlined}):
\begin{center}
\begin{tabular}{lllll}
\toprule
$r$ & computed capacity & exact capacity & relative error & time (s) \\
\midrule
$0.5$ & $\underline{1.03065123518701}5$ & $1.030651235187014$ & $8.62 \cdot
10^{-16}$ & $0.1$ \\
$0.7$ & $\underline{1.25247255560197}0$ & $1.252472555601971$ & $1.42 \cdot
10^{-15}$ & $0.1$ \\
$0.9$ & $\underline{1.46569864072979}5$ & $1.465698640729791$ & $2.42 \cdot
10^{-15}$ & $0.1$ \\
\bottomrule
\end{tabular}
\medskip
\end{center}
}
\end{example}

\begin{example}[Two disks with different radii] \label{ex:2disks_uneq}
{\rm
Let $0 < u < v$ and
\begin{equation*}
a = \frac{\sinh(v)}{\sinh(v-u)} \quad \text{and} \quad r =
\frac{\sinh(u)}{\sinh(v-u)}.
\end{equation*}
Then the capacity of $E_{a,r} = D_1(0) \cup D_r(a)$, which is the union of two
disjoint disks, is
\begin{equation*}
c(E_{a,r}) = e^{u^2/v} \sinh(u) \abslr{ \frac{ \theta_2(0; e^{-v}) \theta_3(0;
e^{-v}) \theta_4(0; e^{-v}) }{ \theta_1(iu; e^{-v}) } },
\end{equation*}
where $\theta_1, \theta_2, \theta_3, \theta_4$ are the Jacobi theta functions
(see below).  This formula was brought to our attention by Thomas Ransford, who
derived it from the Green's function of an annulus
in~\cite[Ch.~V]{CourantHilbert1953}.  We are not aware of a publication of this 
result in the literature.
The Jacobi theta functions are~\cite[§21.1]{WhittakerWatson1962}
\begin{equation*}
\begin{split}
\theta_1(z; q) &= 2 q^{1/4} \sum_{n=0}^\infty (-1)^n q^{n(n+1)} \sin((2n+1) z),
\\
\theta_2(z; q) &= 2 q^{1/4} \sum_{n=0}^\infty q^{n(n+1)} \cos((2n+1) z), \\
\theta_3(z; q) &= 1 + 2 \sum_{n=1}^\infty q^{n^2} \cos(2n z), \\
\theta_4(z; q) &= 1 + 2 \sum_{n=1}^\infty (-1)^n q^{n^2} \cos(2n z),
\end{split}
\end{equation*}
and we evaluate them by truncating the series when the absolute value of
the terms become smaller than MATLAB's \texttt{eps}.  We use the computed value
of $c(E_{a,r})$ as the exact capacity.
We apply our numerical methods with $n = 2^8$ nodes per boundary, and obtain
the very accurate results shown in the following table (correct digits are
\underline{underlined}):
\begin{center}
\begin{tabular}{llllr}
\toprule
$(u, v)$ & computed capacity & exact capacity & relative error & time (s) \\
\midrule
$(0.5, 0.7)$ & $\underline{2.99127153954169}6$ & $2.991271539541695$ & $2.96
\cdot 10^{-16}$ & $0.2$ \\
$(0.5, 1.0)$ & $\underline{1.637069166040759}$ & $1.637069166040759$ & $2.72
\cdot 10^{-16}$ & $0.1$ \\
$(0.5, 1.5)$ & $\underline{1.2602091592322}60$ & $1.260209159232259$ & $1.76
\cdot 10^{-16}$ & $0.1$ \\
\bottomrule
\end{tabular}
\medskip
\end{center}
}
\end{example}

\begin{example} \label{ex:Rostand_disks}
{\rm
Let $E$ be the union of a disk and two half-disks, as shown in
Figure~\ref{fig:bw_hd} (right). We are not aware of an analytic formula for
$c(E)$, but it has been shown in~\cite{RansfordRostand2007} that
\begin{equation*}
c(E)\,\in\,[2.1969933,\,2.2003506],
\end{equation*}
and the authors ``best guess'' is $c(E)\approx 2.19699371717$.
We parameterize the circle by $\eta_1(t) = e^{-\i t}$, and the two half-disks
analogously to Example~\ref{ex:half-disk}. With $n = 2^{10}$ nodes per boundary
component our method gives the computed value $c(E)\approx 2.196993710282112$
in about $0.9~\textrm{s}$. For $n = 2^{16}$ our method computes the value
$c(E)\approx 2.196993717171386$ in $36.2~\textrm{s}$, and this value
matches exactly the estimate from~\cite{RansfordRostand2007}.
}
\end{example}

\begin{figure}
\centerline{
\subfigure[$7$
ellipses]{\includegraphics[width=0.5\textwidth]{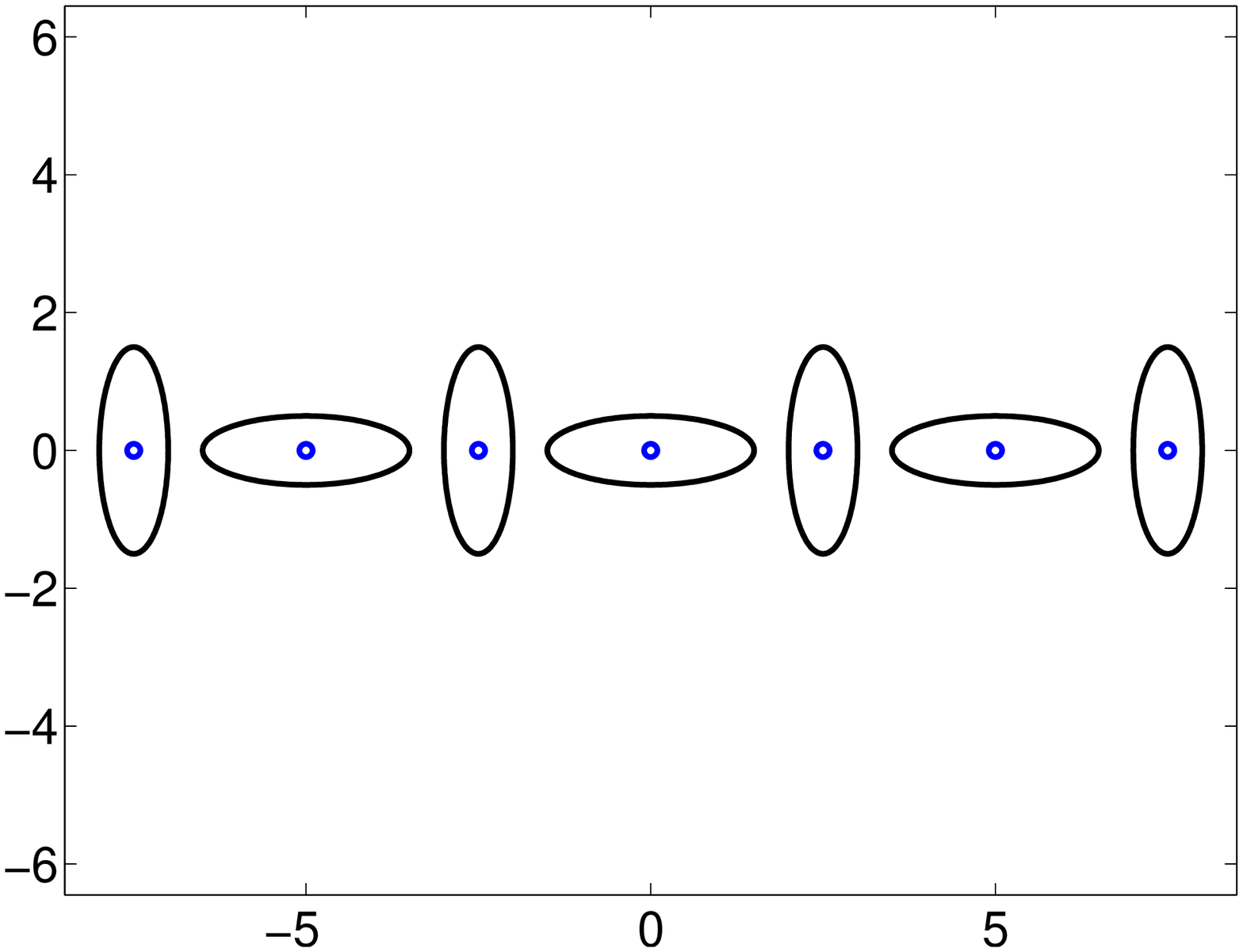}}
\subfigure[$64$ disks]{\includegraphics[width=0.5\textwidth]{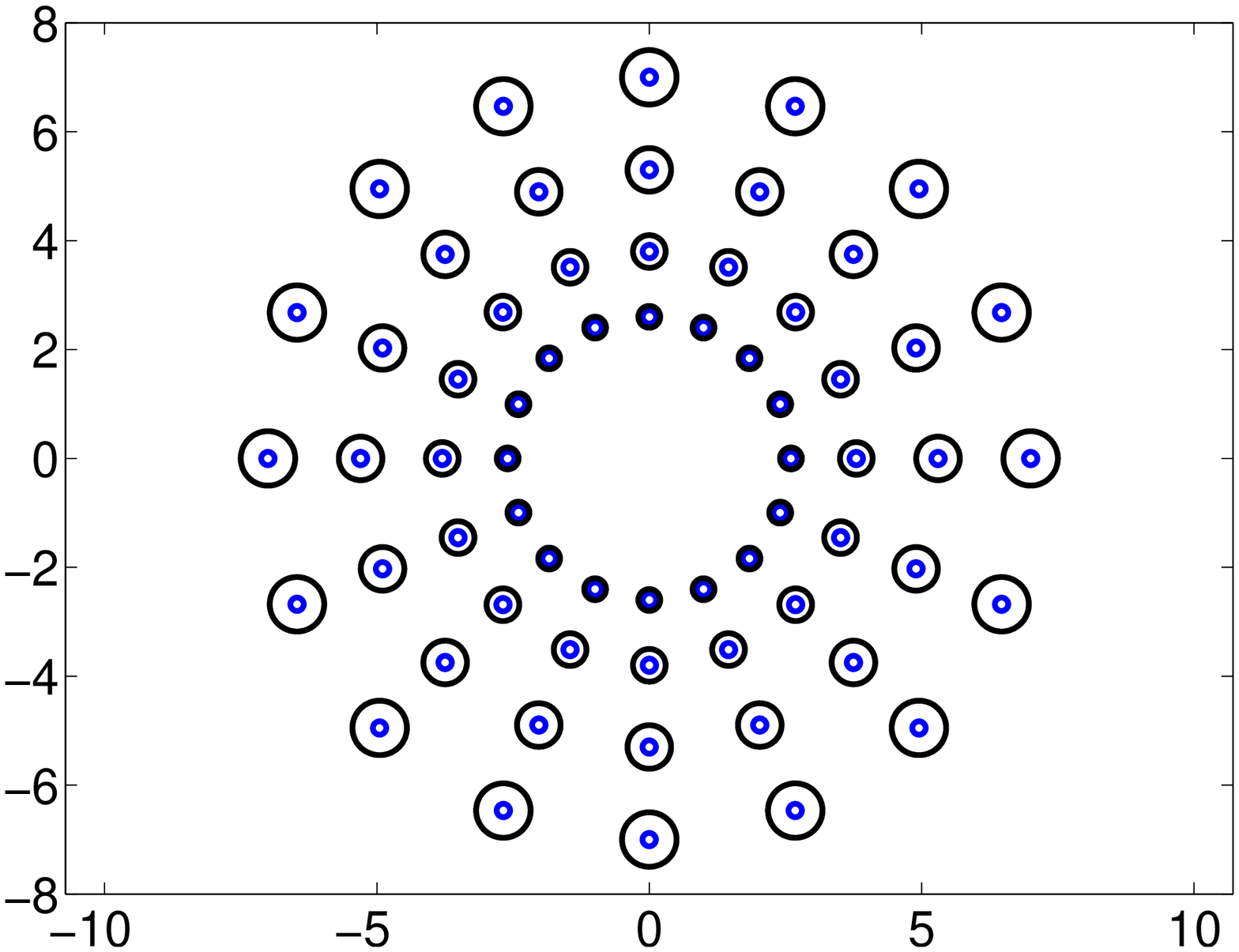}}
}
\centerline{
\subfigure[$4$ squares]{\includegraphics[width=0.5\textwidth]{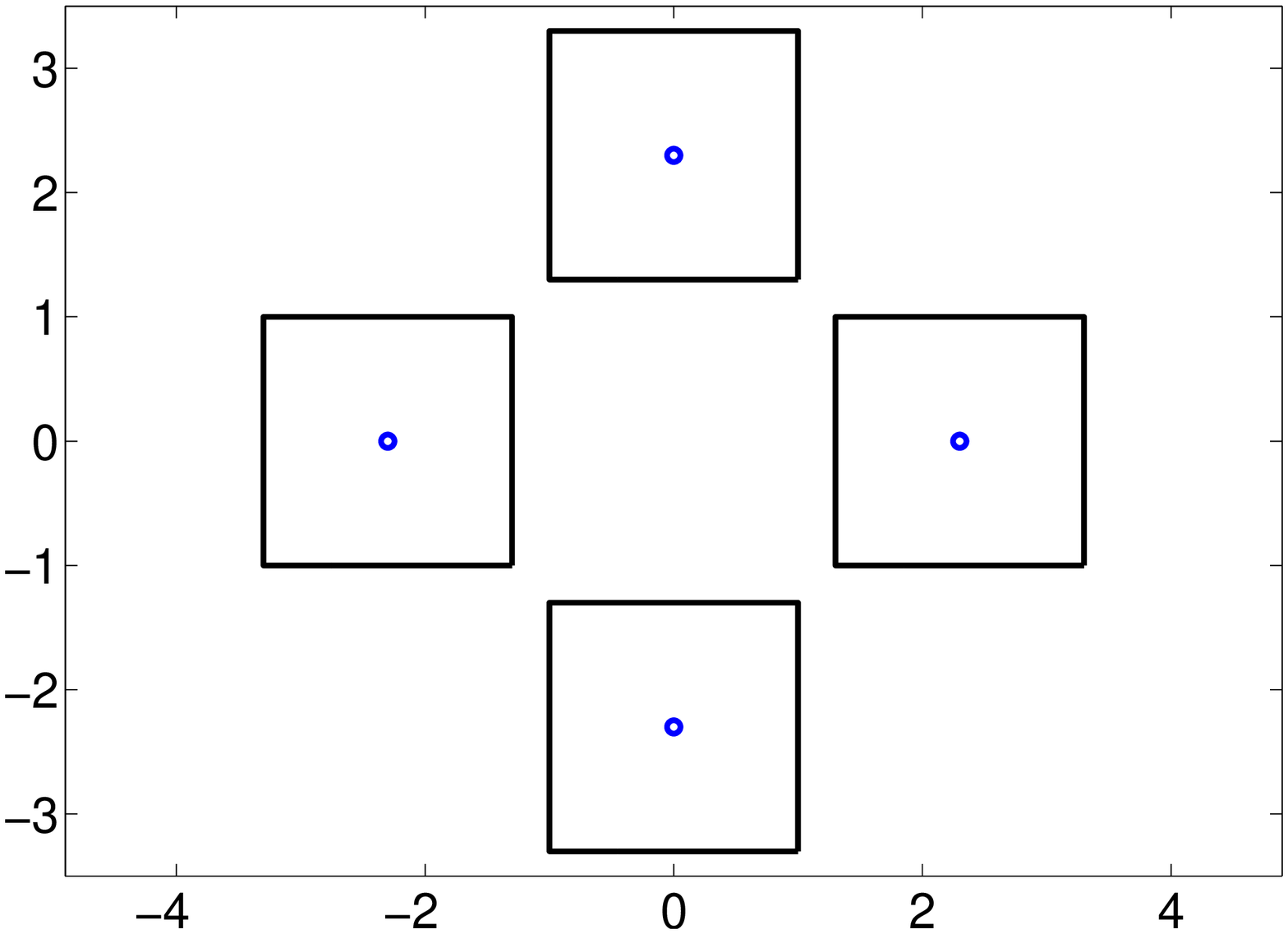}}
\subfigure[{$3$ sets as
in~\cite{KochLiesen2000}}]{\includegraphics[width=0.5\textwidth]{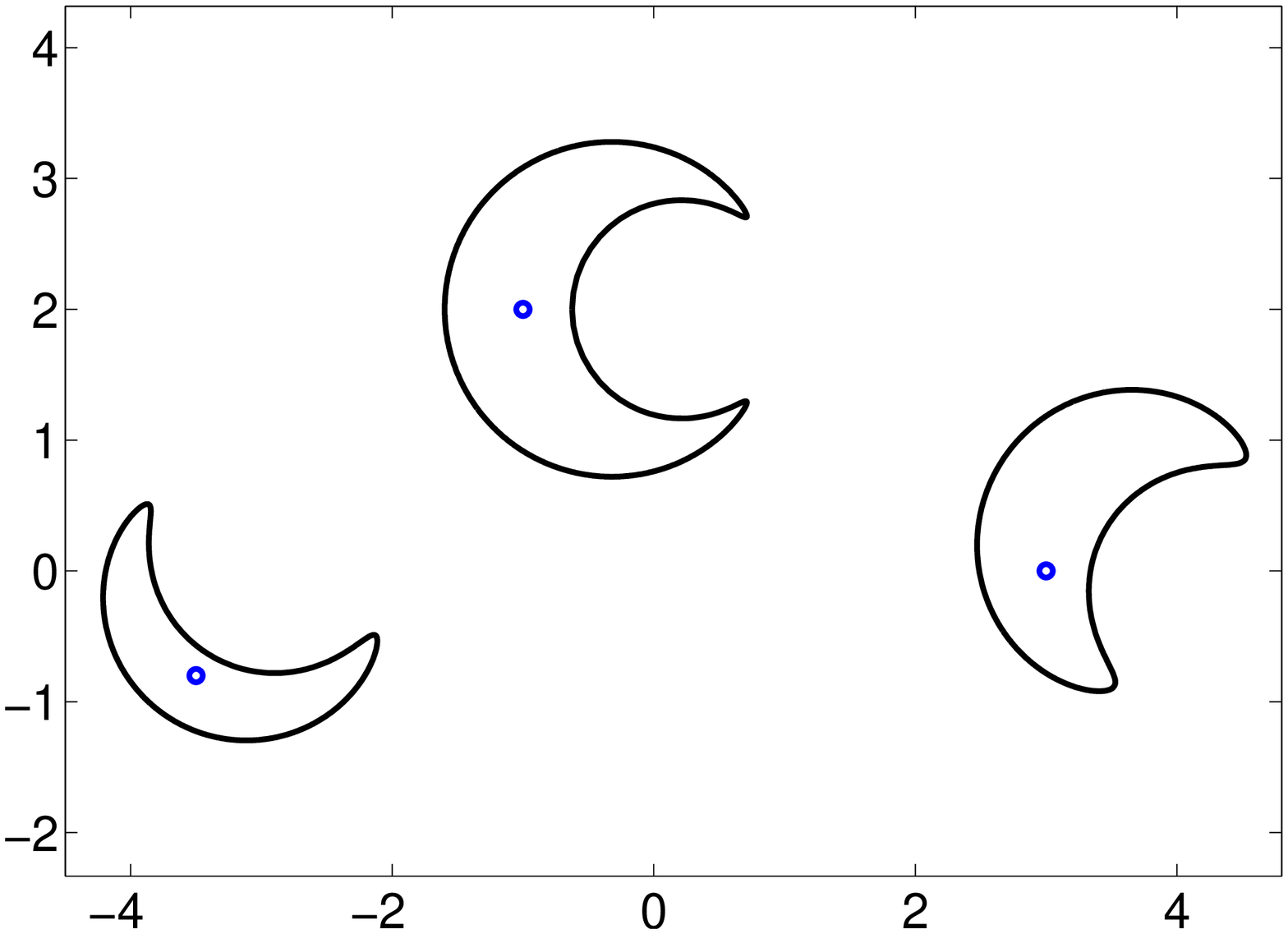}}
}
\caption{The sets in Example~\ref{ex:multcomp}. The (blue) dots show the
auxiliary points $\alpha_j$.}
\label{fig:multcomp}
\end{figure}

\begin{example} \label{ex:multcomp}
{\rm
We consider the compact sets shown in Figure~\ref{fig:multcomp}. These
sets were also used in~\cite[Examples~2--5]{NLS2016_numconf}.  To our
knowledge, the logarithmic capacities of these sets are not known analytically.
The following table shows the values computed by our method and the computation
times:
\begin{center}
\begin{tabular}{lllr}
\toprule
$E$ & $n$ & computed capacity & time (s) \\
\midrule
$7$ ellipses & $2^8$ & $4.961809958325545$ & $1.2$ \\
$64$ disks & $2^8$ & $7.177814562549484$ & $43.6$ \\
$4$ squares & $2^{10}$ & $3.083190170261768$ & $1.3$ \\
$3$ sets as in \cite{KochLiesen2000} & $2^{10}$ & $2.977866214534663$ & $1.1$ \\
\bottomrule
\end{tabular}
\medskip
\end{center}
}
\end{example}

\begin{figure}
\centerline{
\subfigure[An aerial photograph of ``The World Islands'']
{\includegraphics[width=0.5\textwidth]{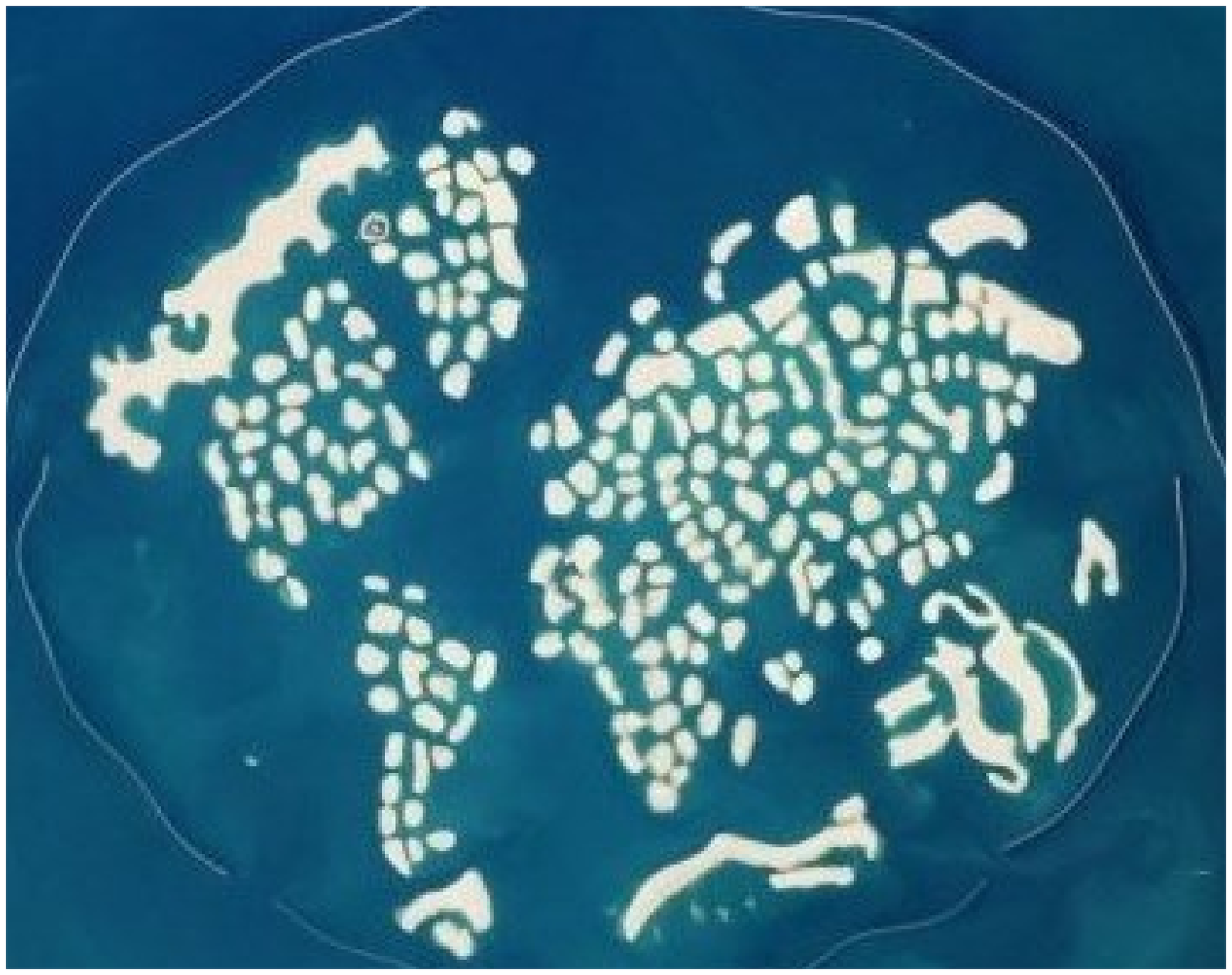}}
\subfigure[The boundaries of the islands extracted from the image]
{\includegraphics[width=0.5\textwidth]{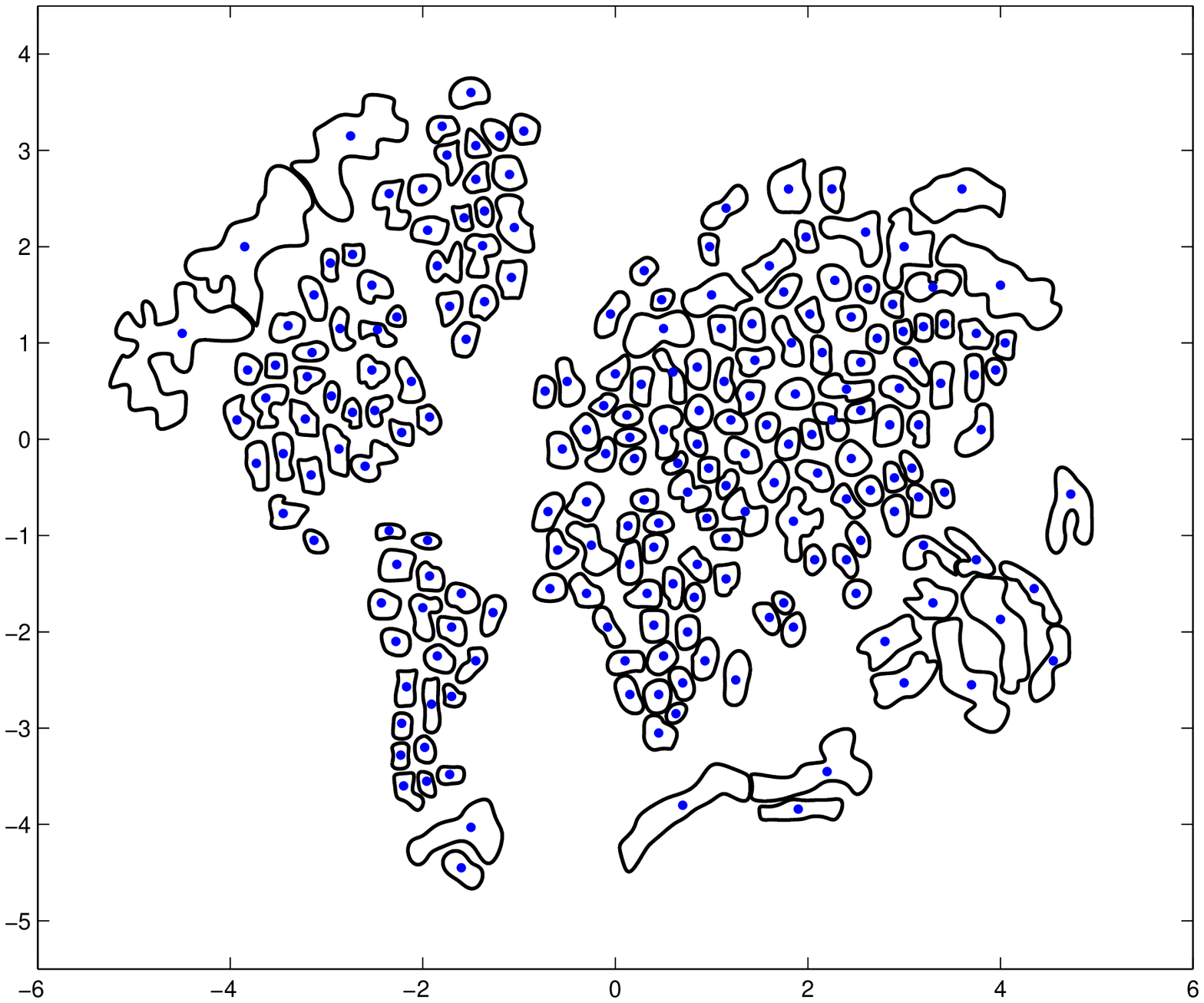}}
}
\caption{``The World Islands'' from Example~\ref{ex:dubai}. The (blue) dots
show the auxiliary points $\alpha_j$.}
\label{fig:dubai}
\end{figure}

\begin{example}[``The World Islands''] \label{ex:dubai}
{\rm
We consider the unbounded domain $\cK$ of connectivity $\ell=210$
exterior to an artificial archipelago located in the waters of the Arabian Gulf,
four kilometres off the coast of Dubai, and known as ``The World Islands''; see
also~\cite{Nas-fast}. An aerial image from {\it Google Maps} of ``The World
Islands''
is shown in Figure~\ref{fig:dubai} (left).  The boundaries of the islands
extracted from the aerial image are shown in Figure~\ref{fig:dubai} (right),
and we parameterize them using trigonometric interpolating polynomials.  The
boundaries are very close to
each other but they do not touch. The following table shows the values computed
by our method and the computation times:
\begin{center}
\begin{tabular}{llr}
\toprule
$n$ & computed capacity & time (s) \\
\midrule
$2^5$    & $4.384226180107323$  & $484$ \\
$2^6$    & $4.388057916386704$  & $872$ \\
$2^7$    & $4.387882300144899$  & $1464$ \\
$2^8$    & $4.387881092385658$  & $2813$ \\
$2^9$    & $4.387881092740317$  & $4986$ \\
$2^{10}$ & $4.387881092740335$  & $9656$ \\
$2^{11}$ & $4.387881092740328$  & $17717$ \\
\bottomrule
\end{tabular}
\medskip
\end{center}
}
\end{example}

\subsection{Several real intervals and Cantor sets}
\label{sect:cantor}

In this section we consider sets consisting of several real intervals, and in 
particular the classical Cantor middle third set and generalizations of it.
Our method is not directly applicable to such sets, since these are not bounded 
by Jordan curves.  We overcome this difficulty by considering a preliminary 
conformal map to ``open up'' the intervals and obtain a compact set of same 
logarithmic capacity, but bounded by smooth Jordan curves.

Let $E$ be a set that consists of $\ell$ real intervals, so that its
complement $\cK = E^c = \widehat{\C} \backslash E$ is an unbounded domain
(containing $\infty$) bounded by $\ell$ parallel straight slits, all on the
real axis.  In order to apply our method, we will first use conformal mappings
to compute
an unbounded multiply connected domain $G$ exterior to a disjoint union
of $\ell$ ellipses, such that the domains $G$ and $\cK$ are conformally
equivalent.  The computation of $G$ and the conformal map $z = \omega(\zeta)$
from $G$ onto $\cK$ is based on a technique developed recently
in~\cite{Nasser-Chris2016}, which we describe in Appendix~\ref{sect:appendix}.
The conformal map is normalized by $\omega(\zeta) = \zeta + O(1/\zeta)$ as 
$\zeta \to
\infty$, which makes it unique, and also implies $c(E) = c(G^c)$
by~\cite[Theorem~5.2.3]{Ransford1995}, that is, the capacity of $E$ can be
computed as the capacity of the union of the $\ell$ ellipses.
In summary, we use the method described in Appendix~\ref{sect:appendix} for
computing $G$, and then apply our usual method (as stated in
Fig.~\ref{fig:logcapacity}) in order to compute an approximation of
$c(E)=c(G^c)$.

We start with the much studied case of two real intervals, for which the
logarithmic capacity is known analytically.  The numerical results for these
sets also give an indication of the accuracy of our method for the Cantor sets,
where no analytic formula for their logarithmic capacity is known.

\begin{example}[Two real intervals]{\rm
Let $E_{a,b}=[-1,a]\cup [b,1]$ with $-1<a<b<1$.  When the two intervals have
the same length, i.e., when $0<-a=b<1$, the exact capacity of $E_{a,b}$ is 
given by (see Table~\ref{table:known_cap})
\[
c(E_{a,b})=\frac{1}{2}\sqrt{1-a^2}.
\]
For the general case, an analytic formula derived by Achieser~\cite{Achieser1930} 
(see also~\cite{Schiefermayr2011}) has the form 
\[
c(E_{a,b})= \frac{1+b}{2(1+a)}\frac{\theta_4(0; q)}{\theta_4(\lambda\pi/2; q)},
\]
where
\[
k=\sqrt{\frac{2(b-a)}{(1-a)(1+b)}}, \quad k'=\sqrt{1-k^2}, \quad K=K(k), \quad 
K'=K(k'), \quad q=e^{-\pi K'/K},
\]
$K(k)$ is the complete elliptic integral of the first kind (see 
Example~\ref{ex:2disks} above), $\theta_4$ is the 
fourth Jacobi theta function (see Example~\ref{ex:2disks_uneq} above), and 
$0<\lambda<K$ is defined uniquely by the Jacobi elliptic function $\sn$ as
\[
\sn(\lambda,k)=\sqrt{\frac{1-a}{2}}.
\]
We compute the values of $K(k)$, $K(k')$, and $\theta_4$ as explained in 
Examples~\ref{ex:2disks} and~\ref{ex:2disks_uneq}.  The value of the parameter 
$\lambda$ is computed using the MATLAB 
function \verb|asne|, i.e.,
\begin{equation*}
\lambda= \verb|asne|(\sqrt{(1-a)/2},k).
\end{equation*}

As explained above, we use our method to compute $c(G^c)=c(E_{a,b})$, where $G$ 
is a domain bounded by smooth Jordan curves that is found by ``opening up'' the 
real intervals.  Very accurate results are obtained with only $n=2^8$ nodes per 
boundary, as shown in the following table (correct digits are 
\underline{underlined}):

\begin{center}
\begin{tabular}{llllr}
\toprule
$(a, b)$ & computed capacity & exact capacity & relative error & time (s) \\
\midrule
$(-0.5, -0.1)$ & $\underline{0.48882927115471}8$ & $0.488829271154715$ & $4.77
\cdot 10^{-15}$ & $2.9$ \\
$(0.5, 0.6)$   & $\underline{0.49910155716636}0$ & $0.499101557166361$ & $1.11
\cdot 10^{-15}$ & $3.5$ \\
$(-0.5, 0.3)$  & $\underline{0.457718411572721}$ & $0.457718411572721$ & $8.49
\cdot 10^{-16}$ & $2.7$ \\
$(-0.5, 0.5)$  & $\underline{0.43301270189221}7$ & $0.433012701892219$ & $5.64
\cdot 10^{-15}$ & $2.1$ \\
$(-0.01, 0.01)$  & $\underline{0.49997499937496}8$ & $0.499974999374969$ & $8.88
\cdot 10^{-16}$ & $3.6$ \\
\bottomrule
\end{tabular}
\medskip
\end{center}

As a final remark concerning this example, it is worth mentioning that an 
analytic formula for the logarithmic capacity of sets consisting of several 
intervals has been derived recently in~\cite{BogatyrevGrigoriev2015}.

}\end{example}

\begin{example}[Cantor middle third set] \label{ex:cantor}
{\rm
In this example we consider the classical Cantor middle third set.  Let
$E_0=[0,1]$ and recursively define
\[
E_k \coloneq \frac{1}{3}E_{k-1}\cup\left(\frac{1}{3}E_{k-1}+\frac{2}{3}\right),
\quad k \geq 1.
\]
This means that $E_k$ is constructed by ``removing'' the middle one third of
each interval that $E_{k-1}$ consists of.  Then the Cantor middle third set is
defined as $E \coloneq \cap_{k=1}^{\infty} E_k$.
While no analytic formula for $c(E)$ is known, several attempts have
been made to numerically approximate $c(E)$. In~\cite{RansfordRostand2007}
it is shown that
\[
c(E)\,\in\, [0.22094810685,\,0.22095089228],
\]
and based on several of their computed values, Ransford and Rostand wrote
that their ``best guess'' is
\begin{equation*}
c(E)\approx 0.220949102189507.
\end{equation*}
Using two different approaches based on Schwarz-Christoffel mappings,
Banjai, Embree, and Trefethen computed $c(E_k)$ for $k=1,2,\dots, 9$,
and extrapolating from their computed values they obtained
$c(E)\approx 0.2209491$\/\footnote{These computations,
made in July 2005, were also reported in~\cite{RansfordRostand2007}.}.
Recently, Kr\"uger and Simon~\cite{KruSim15} obtained the value
$c(E)\approx 0.22094998647421$ in a study of the spectral theory of orthogonal
polynomials associated to the Cantor measure. Referring to their result they
noted that one ``should only trust the first six digits or so''.

We will now describe our approach for computing an approximation of $c(E)$.
Similar to Banjai, Embree and Trefethen, we will compute $c(E_k)$ for a few
small values of $k$, and then obtain an approximation of $c(E)$ by
extrapolation.
We compute the capacities $c(E_K)$ with the open-up method described in the
beginning of Section~\ref{sect:cantor}.  The resulting values and
computation times for $k=1,2,\dots,12$ are shown in the following table:
\begin{center}
\begin{tabular}{rrlr}
\toprule
$k$   &  $\ell=2^k$ & computed capacity    & time (s) \\
\midrule
$1$   &    $2$      & $0.235702260395518$  &     $1.0$ \\
$2$   &    $4$      & $0.228430704425426$  &     $1.5$ \\
$3$   &    $8$      & $0.224752818755436$  &     $2.5$ \\
$4$   &   $16$      & $0.222887290751916$  &     $4.1$ \\
$5$   &   $32$      & $0.221938129124324$  &     $7.3$ \\
$6$   &   $64$      & $0.221454205006181$  &    $15.7$ \\
$7$   &  $128$      & $0.221207178734289$  &    $44.6$ \\
$8$   &  $256$      & $0.221080995391656$  &   $148.0$ \\
$9$   &  $512$      & $0.221016516406108$  &   $565.1$ \\
$10$  & $1024$      & $0.220983561713855$  &  $2375.9$ \\
$11$  & $2048$      & $0.220966717159289$  &  $9128.4$ \\
$12$  & $4096$      & $0.220958106742622$  & $34984.7$ \\
\bottomrule
\end{tabular}
\medskip
\end{center}

In order to extrapolate from our computed values, we note that the differences
\begin{equation*}
d_k=c(E_k)-c(E_{k+1}),\quad k=1,2,\dots,11
\end{equation*}
behave linearly on a logarithmic scale; see the (blue)
circles in Figure~\ref{fig:extrapolation}.
We therefore use the MATLAB command {\tt p=polyfit(1:11,log(d(1:11),1))} for
computing
a linear interpolant $p(x)=p_1 x+p_2$ of the values $\log(d_k)$. The computed
coefficients are
\begin{equation*}
p_1=-0.673356333942526,\quad p_2=-4.26116079806122,
\end{equation*}
and the values $p(k)$ for all $k=1,2,\dots,48$, where $p(48)\approx 10^{-16}$,
are shown by the (black) pluses in Figure~\ref{fig:extrapolation}. Since
$p(k)\approx \log(d_k)$,
we can find an approximation of $c(E_k)$ for each $k=13,14,\dots$ by
extrapolation starting with our computed value for $c(E_{12})$, and obtain
\begin{equation*}
c(E_k) = c(E_{12}) - \sum_{j=12}^{k-1} \exp(p(j)), \quad k \geq 13.
\end{equation*}
Since $p(49)<10^{-16}$, we use only the values up to $p(48)$, which gives our
estimate for the capacity of the Cantor middle third set as
\begin{equation*}
c(E)\approx 0.220949194629475.
\end{equation*}
This estimate agrees up to the seventh digit with the estimates of Ransford and
Rostand as well as Banjai, Embree and Trefethen.

\begin{figure}
\centerline{
\includegraphics[width=0.48\linewidth]{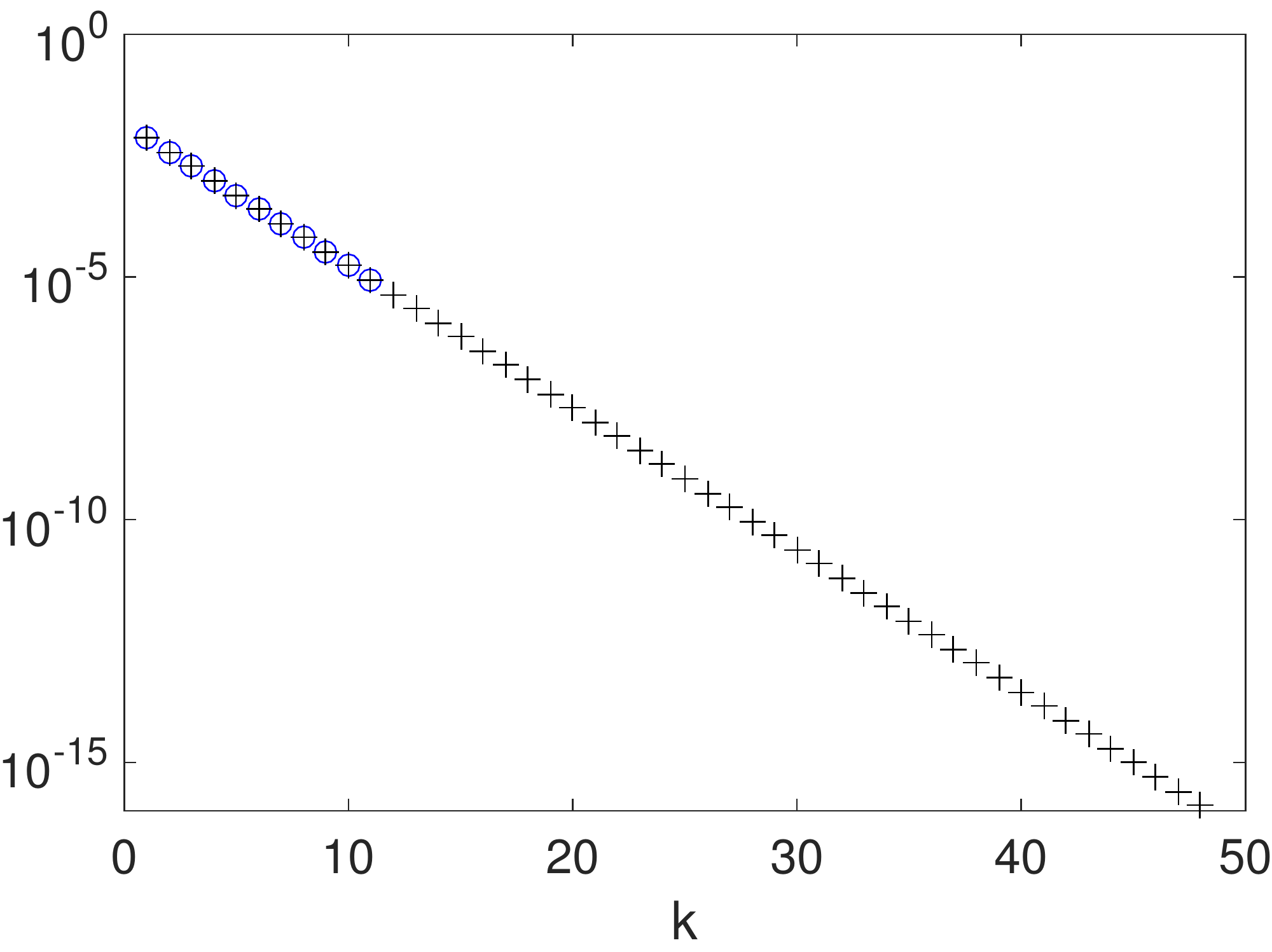}}
\caption{The values $d_k=c(E_{k})-c(E_{k+1})$ for $k=1,2,\dots,11$ (circles)
and the values $\exp(p(k))$ for $k=1,2,\dots,48$ (pluses); see 
Example~\ref{ex:cantor}.}
\label{fig:extrapolation}
\end{figure}

}\end{example}

\begin{figure}
\begin{center}
\includegraphics[width=0.48\linewidth]{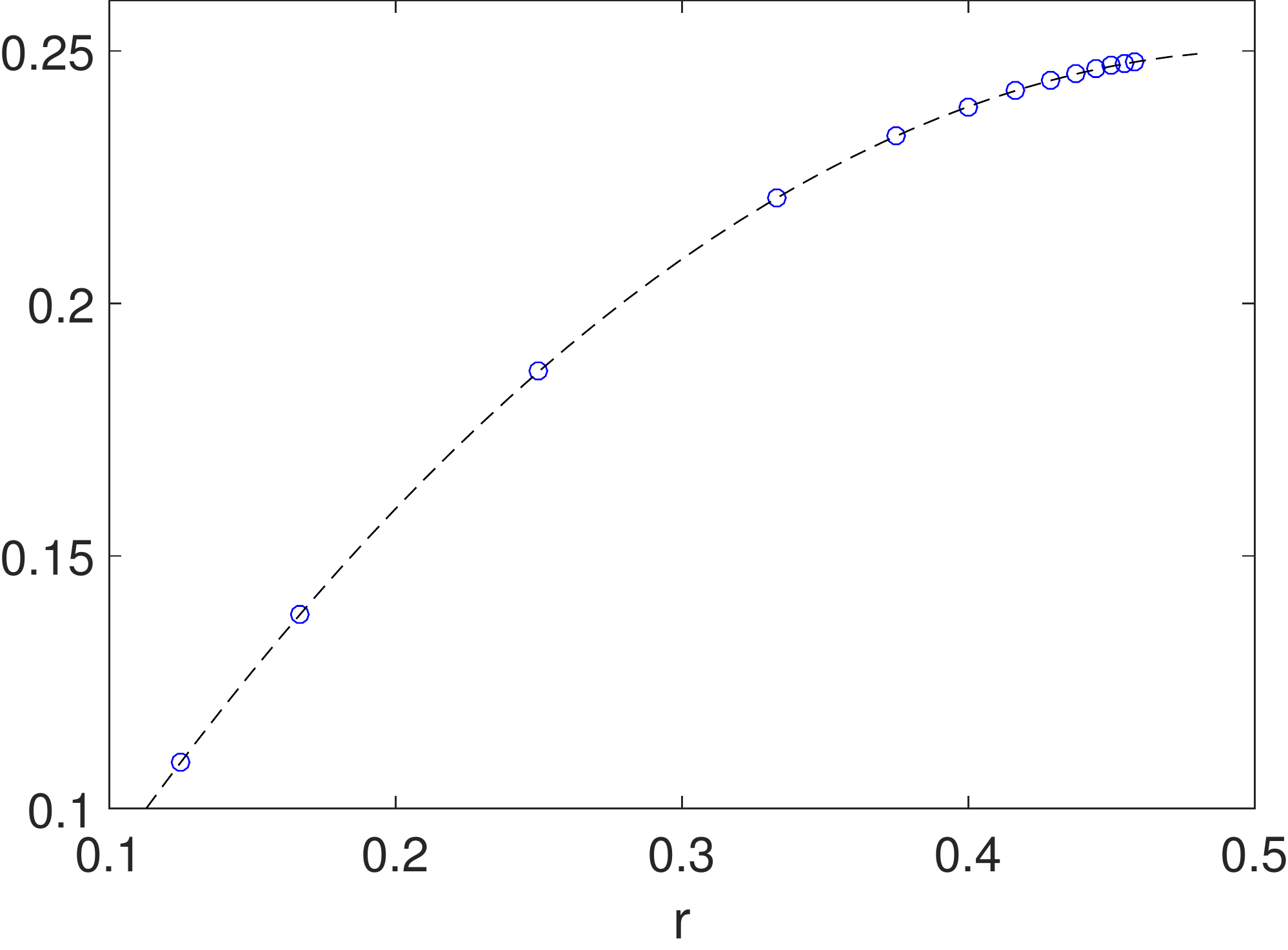}
\includegraphics[width=0.48\linewidth]{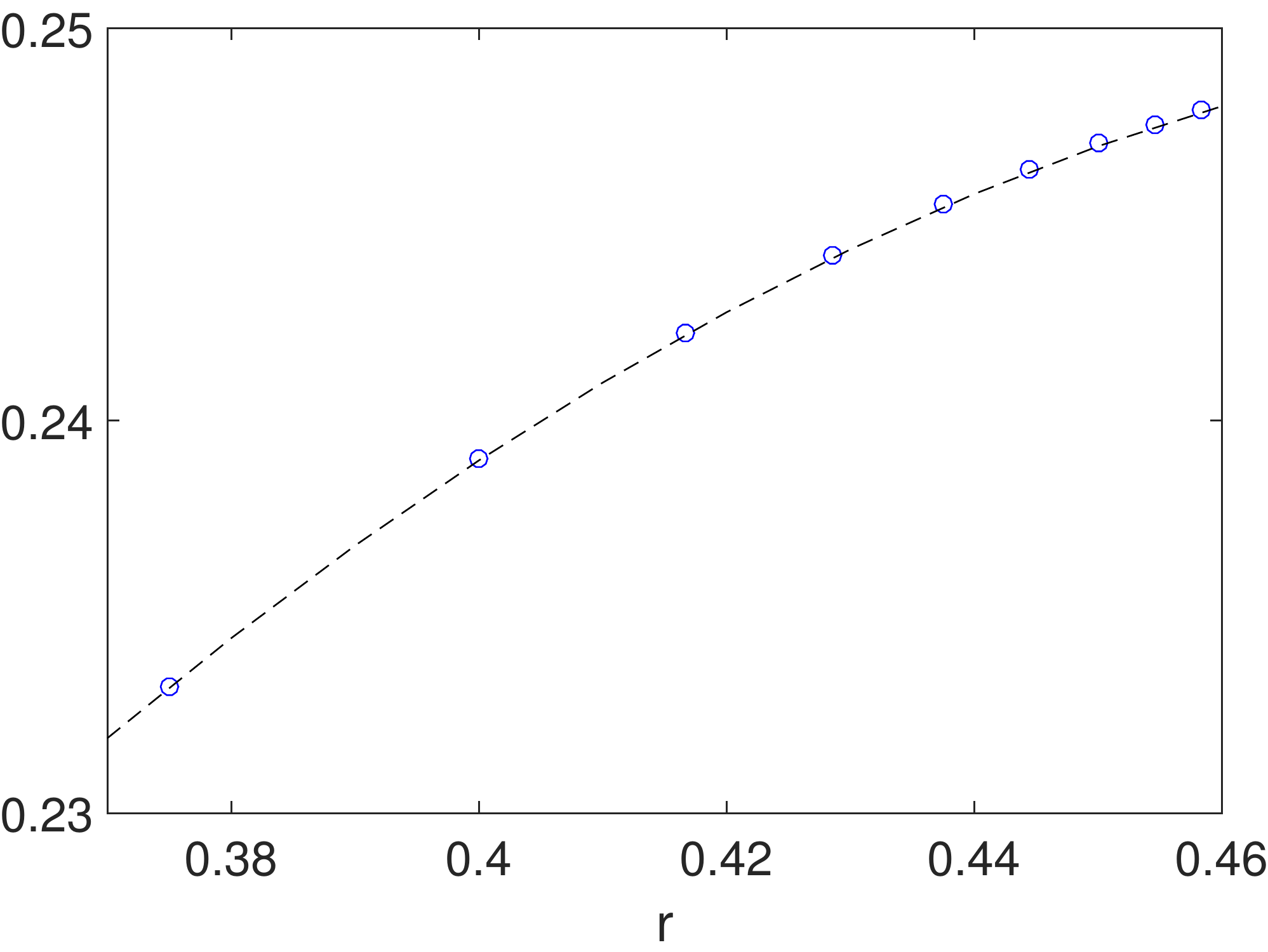}
\end{center}
\caption{Computed approximation of $c(E^r)$ (circles) and the function $f(r)$
(dashed); see Example~\ref{ex:Cantor_gen}.}
\label{fig:cantor_general}
\end{figure}

\begin{example}[Generalized Cantor set]\label{ex:Cantor_gen}
{\rm
As in~\cite[Section~6]{RansfordRostand2007} we will now generalize the
construction of the Cantor middle third set as follows. Let $r\in(0,\,0.5)$ and
$E_0^r \coloneq [0,1]$.
Recursively define
\begin{equation*}
E_k^r \coloneq r E_{k-1}^r \cup \left(r E_{k-1}^r + (1-r)\right), \quad k \geq
1,
\end{equation*}
and let $E^r\coloneq \cap_{k=0}^\infty E_k^r$.  The parameter $r$ determines
how much is removed from each interval of $E_{k-1}^r$ in order to obtain
$E_k^r$.  If we want to remove the middle $q\in (0,1)$, we need to set
$r=(1-q)/2$.  For $q=1/3$ we have $r=1/3$ and hence $E^{1/3}$ is the classical
Cantor middle third set.  The limiting cases are
$E^0=\{0,1\}$ with $c(E^0)=0$, and $E^{1/2}=[0,1]$ with $c(E^{1/2})=1/4$; see
Table~\ref{table:known_cap}.

Using exactly the same approach as described in Example~\ref{ex:cantor} we have
computed the following approximations of $c(E^r)$:

\begin{center}
\begin{tabular}{ccl}
\toprule
$q$ & $r$   &  computed capacity   \\
\hline
$3/4$ & $1/8$ & $0.109156838696175$ \\
$2/3$ & $1/6$ & $0.13844418298159$   \\
$1/2$ & $1/4$ & $0.186511016338442$  \\
$1/3$ & $1/3$ & $0.220949194629475$  \\
$1/4$ & $3/8$ & $0.233218551525021$  \\
$1/5$ & $2/5$ & $0.23901897053678$   \\
$1/6$ & $5/12$ & $0.242233234580321$  \\
$1/7$ & $3/7$ & $0.244206003640726$  \\
$1/8$ & $7/16$ & $0.245506481568117$  \\
$1/9$ & $4/9$ & $0.246410328817$  \\
$1/10$ & $9/20$ & $0.247064652445187$  \\
$1/11$ & $10/22$ & $0.247553947239903$  \\
$1/12$ & $11/24$ & $0.247929630663845$  \\
\bottomrule
\end{tabular}
\medskip
\end{center}

The (blue) circles in Figure~\ref{fig:cantor_general} show our computed
approximations of $c(E^r)$, where the right part of the figure is a
closeup of the left part. The dashed line shows the function
\begin{equation*}
f(r)=r(1-r)-\frac{r^3}{2}\left(\frac12-r\right)^{3/2},
\end{equation*}
which was suggested in~\cite{RansfordRostand2007} as an approximation
of $c(E^r)$. The maximum distance between the values of $f(r)$ and our
computed approximations of $c(E^r)$ is $7.5189 \cdot 10^{-5}$.  Thus,
similar to computations reported in~\cite{RansfordRostand2007}, the function
$f(r)$ very closely approximates our computed approximations of $c(E^r)$.
}\end{example}

\section{Concluding remarks}

We have presented a numerical method for the computation of the logarithmic
capacity of compact sets bounded by Jordan curves in the complex plane.  These
sets may consist of several components and need not have any special symmetry
properties.  In several numerical examples with sets for which the logarithmic
capacity is known analytically, our method yields a computed approximation with
a relative error close to the machine precision.
For ``simple'' sets, in particular simply connected ones, the computations 
in MATLAB take at most a few seconds.

Let us point out a few open questions.  From a computational point of view, an 
automated choice of the auxiliary points $\alpha_j$ (interior to each boundary 
curve) would be of interest.
The numerical experiments for compact sets where the logarithmic capacity is
known analytically suggest that the method is fast and accurate.  A formal
analysis of the numerical stability and accuracy of our method is beyond the
scope of this article and remains a subject of further work.
To compute an approximation of the capacity
of the Cantor sets, we devised an ad hoc method to ``open up'' the intervals by
conformal mapping and obtain a domain bounded by Jordan curves, to which our
method could be applied.  It would be of interest to devise an ``open-up 
method'' for general Jordan arcs that is computationally tractable.

\appendix
\section{Numerical computation of the preimage of a parallel slit domain}
\label{sect:appendix}

\begin{figure}[t] \centerline{
\scalebox{0.5}{\includegraphics{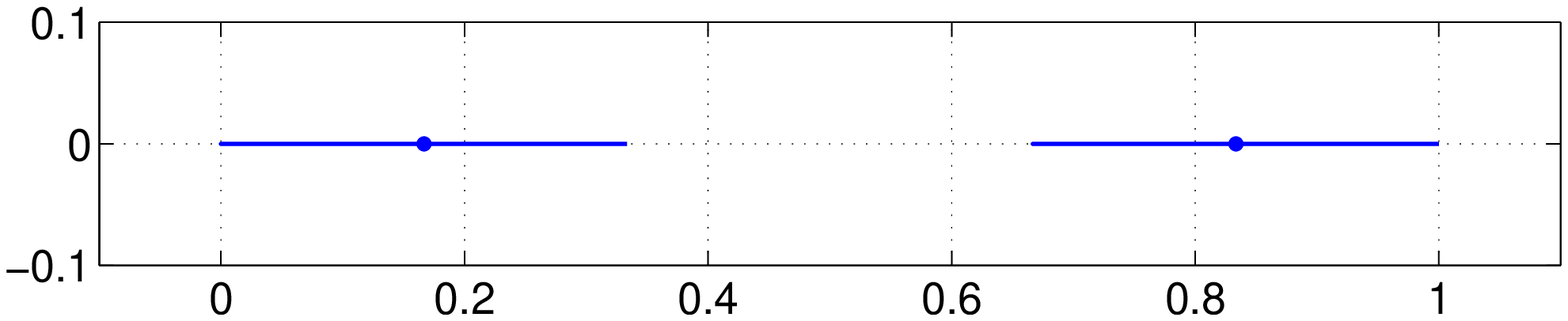}}
}
\centerline{
\scalebox{0.5}{\includegraphics{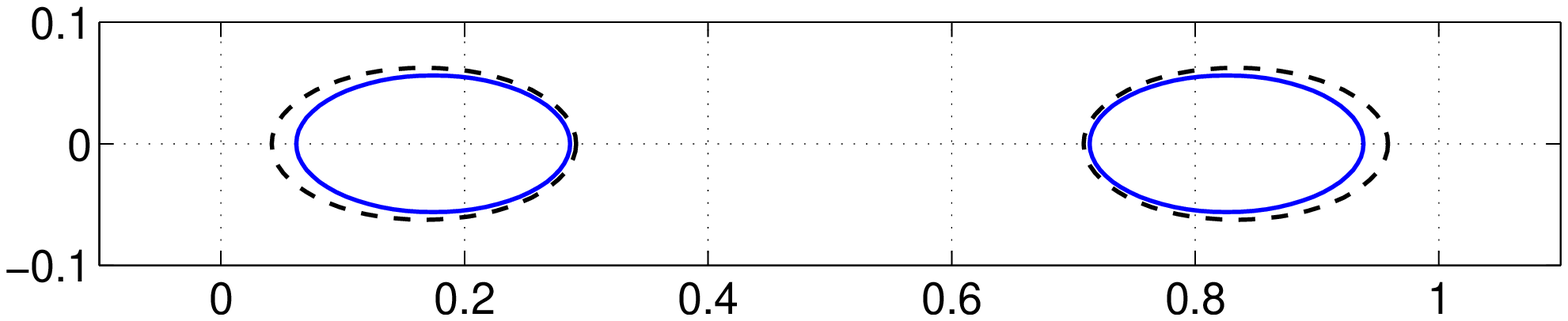}}
}
\caption{The given parallel slit domain $\Omega$ (top), the initial preimage
domain $G^0$ (dashed line, bottom), and the computed preimage domain $G$
(solid line, bottom).}
\label{f:str}
\end{figure}

Let $\Omega$ be a given parallel slit domain, i.e., the entire $z$-plane
with $m$ slits $L_j$, $j=1,2,\ldots,\ell$, along straight lines; see the top
of Figure~\ref{f:str}. An efficient numerical method for computing the
conformal map $z=\omega(\zeta)$ from an unbounded domain $G$ exterior
to $\ell$ smooth Jordan curves $\Gamma_j$, $j=1,2,\ldots,\ell$, onto the
parallel slit domain $\Omega$ such that
$\omega(\zeta) = \zeta + O(1/\zeta)$ as $\zeta \to \infty$
has been presented in~\cite[Section~4.5]{Nasser-SIAM2009}. Assume that the
boundary $\Gamma$ of $G$ is parametrized by the function $\eta(t)$ as
in~\eqref{eqn:eta}. Assume also the operators $\bN$ and $\bM$ are the
same operators as in~\eqref{eqn:bieq}.  Then we have the following theorem
from~\cite{Nasser-SIAM2009}.

\begin{theorem}
\label{thm:cm-app}
Let
\begin{equation}
\gamma(t)={\mathrm{Im}}[\eta(t)], \quad t\in J,
\label{eqn:gama-app}
\end{equation}
let $\mu$ be the unique solution of the boundary integral equation
\begin{equation}
(\bI - \bN) \mu = - \bM \gamma,
\label{eqn:bieq-app}
\end{equation}
and let $h$ be the piecewise constant function
\begin{equation}
h = ( \bM \mu - (\bI - \bN) \gamma)/2.
\label{eqn:h-app}
\end{equation}
Then the function $f$ with the boundary values
\begin{equation}
f(\eta(t))=\gamma(t)+h(t)+\i\mu(t)
\label{eqn:f-app}
\end{equation}
is analytic in $G$ with $f(\infty)=0$ and the conformal mapping $\omega$ is 
given by
\begin{equation}
\omega(\zeta)=\zeta-\i f(\zeta), \quad \zeta\in G\cup\Gamma.
\label{eqn:omega-app}
\end{equation}
\end{theorem}

In Theorem~\ref{thm:cm-app}, the domain $G$ is assumed to be known,
and the integral equation~\eqref{eqn:bieq-app} is used to the find
the conformal map $z=\omega(\zeta)$ from $G$ onto the parallel slit domain
$\Omega=\omega(G)$. In our application with the Cantor sets, however,
the domain $G$ is unknown and the parallel slit domain $\Omega$ is known.
Hence a straightforward application of a numerical method based on
Theorem~\ref{thm:cm-app} is not possible.

We will now describe an iterative method developed in~\cite{Nasser-Chris2016}
for computing $G$ and the conformal map from $G$ onto the (known) parallel
slit domain $\Omega$.  The method is an improvement of a numerical method
suggested by Aoyama, Sakajo, and Tanaka~\cite{Aoy-Sak-Tan13}, where
the preimage $G$ is assumed to be circular.  Since the image region $\Omega$
is elongated (parallel slit domain), crowding can cause serious problems.
Further, the convergence of the iterative method is slow if $G$ is
assumed to be circular. To overcome such difficulties, it was assumed
in~\cite{Nasser-Chris2016} that the boundaries of the domain $G$ are
ellipses instead of circles.

Let $|L_j|$ denote the length of the slit $L_j$ and let $z_j$ denote its center,
$j=1,2,\ldots,\ell$. In the iteration step $i=0,1,2,\ldots$ we assume that the
domain $G^i$ is a multiply connected domain bounded by the ellipses $\Gamma^i_j$,
parametrized for $j=1,2,\ldots,\ell$ by
\[
\eta^i_j(t)=\zeta^i_j+0.5(a^i_j\cos t-\i b^i_j\sin t), \quad t\in J_j=[0,2\pi].
\]
Then the following iteration computes the centers of the ellipses $\zeta^i_j$,
the lengths of the major axes $a^i_j$, and the lengths of the minor axes $b^i_j$
for $j=1,2,\ldots,\ell$.

\medskip

\noindent{\bf  Initialization:}\\
Let $\varepsilon>0$ be a given tolerance and let $Max$ be a maximum number of 
iterations.
(In our numerical experiments in this paper we always used
$\varepsilon=10^{-14}$ and $Max=50$.) Set
\[
\zeta^0_j=z_j, \quad a^0_j=(1-0.5r)|L_j|, \quad b^0_j=r a^0_j,
\]
where $0<r<1$ is the ratio of the lengths of the major and minor axes of the
ellipse (see Figure~\ref{f:str}
(dashed line, bottom) for $r=0.5$).

\noindent{\bf For $i=1,2,\dots$:}
\begin{enumerate}
\item Map $G^{i-1}$ to a parallel slit domain $\Omega^i$ (based on 
Theorem~\ref{thm:cm-app}),
which is the entire $z-$plane with $\ell$ slits $L^i_j$, $j=1,2,\ldots,\ell$,
along horizontal straight lines.
\item If $|L^i_j|$ denotes the length of the slit $L^i_j$
and $z^i_j$ denotes its center, then we define the parameters
of the preimage domain $G^i$ as
\begin{eqnarray}
\label{eq:slt-k}
\zeta^{i}_j &=& \zeta^{i-1}_j-(z^{i}_j-z_j), \\
a^{i}_j &=& a^{i-1}_j-(|L^{i}_j| -|L_j|), \\
b^{i}_j &=& r a^{i}_j.
\end{eqnarray}
\item Stop the iteration if
\[
\max_{1\le j\le m}(|z^{i}_j-z_j|+||L^{i}_j| -|L_j||)<\varepsilon \quad{\rm or} 
\quad i>Max
\]
\end{enumerate}

Several numerical examples in this paper as well as 
in~\cite{Aoy-Sak-Tan13,Nasser-Chris2016} show the convergence of this iterative 
method, but no proof of convergence has been given so far. Numerical 
experiments also show that the iterative method requires fewer iterations
for small values of $r$, which means that the ellipses will be thin. For thin 
ellipses, however, we usually need a larger number of points $n$ for 
discretizing the boundary integral equations and the GMRES method for solving 
these discretized equations requires more iterations to converge.

\begin{figure}
\centerline{
\includegraphics[width=0.48\linewidth]{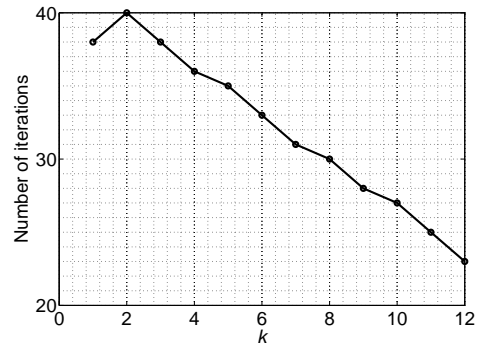}}
\caption{Number of iterations for computing the preimage domains required in 
the computation of $c(E_k)$ for $k=1,2,\dots,12$.}
\label{f:itr}
\end{figure}

In the numerical experiments with the Cantor sets shown in this paper we have 
not chosen to optimize upon these parameters, but we used the fixed values 
$r=0.5$ and $n=64$.  The number of iterations for the convergence to the 
accuracy $\varepsilon=10^{-14}$ of the above iterative method applied in the 
computation of $c(E_k)$ for $k=1,2,\dots,12$ is shown in Figure~\ref{f:itr}.
The (unpreconditioned) GMRES method for solving the discretized
integral equations required between $5$ and $11$ iterations.

\paragraph{Acknowledgments}
We thank Nick Trefethen for sharing the numerical results on the capacity of 
the Cantor middle third set he obtained together with Banjai and Embree.
We thank Thomas Ransford for bringing to our attention the analytic formula for 
the capacity of two unequal disks (Example~\ref{ex:2disks_uneq}).

\bibliographystyle{siam}
\bibliography{capacity}

\begin{thebibliography}{10}

\bibitem{Achieser1930}
{\sc N.~{Achieser}}, {\em {Sur les polynomes de Tchebycheff pour deux
  segments.}}, {C. R. Acad. Sci., Paris}, 191 (1930), pp.~754--756.

\bibitem{Aoy-Sak-Tan13}
{\sc N.~Aoyama, T.~Sakajo, and H.~Tanaka}, {\em A computational theory for
  spiral point vortices in multiply connected domains with slit boundaries},
  Jpn. J. Ind. Appl. Math., 30 (2013), pp.~485--509.

\bibitem{BogatyrevGrigoriev2015}
{\sc A.~B. {Bogatyrev} and O.~A. {Grigoriev}}, {\em Capacity of several aligned
  segments}, ArXiv: 1512.07154,  (2015).

\bibitem{CourantHilbert1953}
{\sc R.~Courant and D.~Hilbert}, {\em Methods of mathematical physics. {V}ol.
  {I}}, Interscience Publishers, Inc., New York, N.Y., 1953.

\bibitem{Davis1957}
{\sc P.~{Davis}}, {\em {Numerical computation of the transfinite diameter of
  two collinear line segments.}}, {J. Res. Natl. Bur. Stand.}, 58 (1957),
  pp.~155--156.

\bibitem{Dijkstra2009}
{\sc W.~Dijkstra and M.~Hochstenbach}, {\em Numerical approximation of the
  logarithmic capacity}, CASA report,  (2009), pp.~08--09.

\bibitem{SCtoolbox}
{\sc T.~A. Driscoll}, {\em Schwarz--{C}hristoffel {T}oolbox {U}ser's {G}uide.
  {V}ersion 2.3}, available from
  \texttt{http://www.math.udel.edu/$\sim$driscoll/SC/guide.pdf}.

\bibitem{DriscollTrefethen2002}
{\sc T.~A. Driscoll and L.~N. Trefethen}, {\em Schwarz-{C}hristoffel mapping},
  vol.~8 of Cambridge Monographs on Applied and Computational Mathematics,
  Cambridge University Press, Cambridge, 2002.

\bibitem{EmbTre99}
{\sc M.~Embree and L.~N. Trefethen}, {\em Green's functions for multiply
  connected domains via conformal mapping}, SIAM Rev., 41 (1999), pp.~745--761.

\bibitem{Fek23}
{\sc M.~Fekete}, {\em \"{U}ber die {V}erteilung der {W}urzeln bei gewissen
  algebraischen {G}leichungen mit ganzzahligen {K}oeffizienten}, Math. Z., 17
  (1923), pp.~228--249.

\bibitem{Gre-Gim12}
{\sc L.~Greengard and Z.~Gimbutas}, {\em {FMMLIB2D}: A {MATLAB} toolbox for
  fast multipole method in two dimensions, Version 1.2,
  \texttt{http://www.cims.nyu.edu/cmcl/fmm2dlib/fmm2dlib.html}}, 2012.

\bibitem{KochLiesen2000}
{\sc T.~Koch and J.~Liesen}, {\em The conformal ``bratwurst'' maps and
  associated {F}aber polynomials}, Numer. Math., 86 (2000), pp.~173--191.

\bibitem{Kress1990}
{\sc R.~Kress}, {\em A {N}ystr\"om method for boundary integral equations in
  domains with corners}, Numer. Math., 58 (1990), pp.~145--161.

\bibitem{Kress1991}
\leavevmode\vrule height 2pt depth -1.6pt width 23pt, {\em Boundary integral
  equations in time-harmonic acoustic scattering}, Math. Comput. Modelling, 15
  (1991), pp.~229--243.

\bibitem{KruSim15}
{\sc H.~Kr{\"u}ger and B.~Simon}, {\em Cantor polynomials and some related
  classes of {OPRL}}, J. Approx. Theory, 191 (2015), pp.~71--93.

\bibitem{Landkof1972}
{\sc N.~S. Landkof}, {\em Foundations of modern potential theory},
  Springer-Verlag, New York-Heidelberg, 1972.
\newblock Translated from the Russian by A. P. Doohovskoy, Die Grundlehren der
  mathematischen Wissenschaften, Band 180.

\bibitem{Nasser-SIAM2009}
{\sc M.~M.~S. Nasser}, {\em Numerical conformal mapping via a boundary integral
  equation with the generalized neumann kernel}, SIAM J. Sci. Comput., 31
  (2009), pp.~1695--1715.

\bibitem{Nas-fast}
\leavevmode\vrule height 2pt depth -1.6pt width 23pt, {\em Fast solution of
  boundary integral equations with the generalized {N}eumann kernel}, Electron.
  Trans. Numer. Anal., 44 (2015), pp.~189--229.

\bibitem{Nas-siam13}
{\sc M.~M.~S. Nasser and F.~A.~A. Al-Shihri}, {\em A fast boundary integral
  equation method for conformal mapping of multiply connected regions}, SIAM J.
  Sci. Comput., 35 (2013), pp.~A1736--A1760.

\bibitem{Nasser-Chris2016}
{\sc M.~M.~S. Nasser and C.~C. Green}, {\em Fast numerical methods for
  describing ideal fluid flow in domains with multiple stirrers}, in
  preparation,  (2016).

\bibitem{NLS2016_numconf}
{\sc M.~M.~S. Nasser, J.~Liesen, and O.~S{\`e}te}, {\em {Numerical computation
  of the conformal map onto lemniscatic domains}}, Comput. Methods Funct.
  Theory,  (2016), pp.~1--27.

\bibitem{Nas-corner}
{\sc M.~M.~S. Nasser, A.~H.~M. Murid, and Z.~Zamzamir}, {\em A boundary
  integral method for the {R}iemann-{H}ilbert problem in domains with corners},
  Complex Var. Elliptic Equ., 53 (2008), pp.~989--1008.

\bibitem{Ransford1995}
{\sc T.~Ransford}, {\em Potential theory in the complex plane}, vol.~28 of
  London Mathematical Society Student Texts, Cambridge University Press,
  Cambridge, 1995.

\bibitem{Ransford2010}
\leavevmode\vrule height 2pt depth -1.6pt width 23pt, {\em Computation of
  logarithmic capacity}, Comput. Methods Funct. Theory, 10 (2010),
  pp.~555--578.

\bibitem{RansfordRostand2007}
{\sc T.~Ransford and J.~Rostand}, {\em Computation of capacity}, Math. Comp.,
  76 (2007), pp.~1499--1520.

\bibitem{Rathsfeld93}
{\sc A.~Rathsfeld}, {\em Iterative solution of linear systems arising from the
  {N}ystr\"om method for the double-layer potential equation over curves with
  corners}, Math. Methods Appl. Sci., 16 (1993), pp.~443--455.

\bibitem{Rostand1997}
{\sc J.~Rostand}, {\em Computing logarithmic capacity with linear programming},
  Experiment. Math., 6 (1997), pp.~221--238.

\bibitem{Saf10}
{\sc E.~B. Saff}, {\em Logarithmic potential theory with applications to
  approximation theory}, Surv. Approx. Theory, 5 (2010), pp.~165--200.

\bibitem{Schiefermayr2011}
{\sc K.~Schiefermayr}, {\em Estimates for the asymptotic convergence factor of
  two intervals}, J. Comput. Appl. Math., 236 (2011), pp.~28--38.

\bibitem{SeteLiesen2016_conf}
{\sc O.~S{\`e}te and J.~Liesen}, {\em On conformal maps from multiply connected
  domains onto lemniscatic domains}, Electron. Trans. Numer. Anal., 45 (2016),
  pp.~1--15.

\bibitem{SeteLiesen2016_fwprop}
\leavevmode\vrule height 2pt depth -1.6pt width 23pt, {\em {Properties and
  Examples of Faber--Walsh Polynomials}}, Comput. Methods Funct. Theory,
  (2016), pp.~1--27.

\bibitem{Sze24}
{\sc G.~Szeg{\"o}}, {\em Bemerkungen zu einer {A}rbeit von {H}errn {M}.
  {F}ekete: \"{U}ber die {V}erteilung der {W}urzeln bei gewissen algebraischen
  {G}leichungen mit ganzzahligen {K}oeffizienten}, Math. Z., 21 (1924),
  pp.~203--208.

\bibitem{Walsh1956}
{\sc J.~L. Walsh}, {\em On the conformal mapping of multiply connected
  regions}, Trans. Amer. Math. Soc., 82 (1956), pp.~128--146.

\bibitem{WhittakerWatson1962}
{\sc E.~T. Whittaker and G.~N. Watson}, {\em A course of modern analysis. {A}n
  introduction to the general theory of infinite processes and of analytic
  functions: with an account of the principal transcendental functions}, Fourth
  edition. Reprinted, Cambridge University Press, New York, 1962.

\end{thebibliography}

\end{document}